\documentclass[3p]{article}

\usepackage{lipsum}
\usepackage{amsfonts}
\usepackage{graphicx}
\usepackage{epstopdf}
\usepackage{algorithm, algpseudocode}

\usepackage{authblk}

\ifpdf
  \DeclareGraphicsExtensions{.eps,.pdf,.png,.jpg,.svg}
\else
  \DeclareGraphicsExtensions{.eps}
\fi

\newcommand{\Deca}{\texttt{Decapodes.jl} }

\usepackage{amsopn}

\newcommand{\R}{\mathbb{R}}

\DeclareMathOperator{\El}{El}

\newcommand{\cat}[1]{\mathsf{#1}}

\usepackage{todonotes}
\usepackage{caption}
\usepackage{subcaption}
\usepackage{nomencl}
\usepackage{booktabs}

\usepackage{amsmath}

\usepackage{geometry}
\usepackage{fancyhdr}
\usepackage{multicol}
\usepackage{bm}
\usepackage[export]{adjustbox}
\usepackage{tikz}
\usepackage{tikz-cd}
\usepackage{assets/styles/tikzit}

\tikzstyle{variable_node}=[fill={rgb,255: red,108; green,154; blue,195}, draw=black, shape=circle, text=white]
\tikzstyle{block}=[fill={rgb,255: red,226; green,143; blue,65}, draw=black, shape=rectangle, text=white]
\tikzstyle{Ob}=[fill=none, draw=none, shape=rectangle]
\tikzstyle{table}=[fill={rgb,255: red,224; green,229; blue,205}, draw={rgb,255: red,95; green,96; blue,98}, shape=rectangle, rounded corners, minimum height=5ex]
\tikzstyle{title}=[fill=white, draw={rgb,255: red,128; green,190; blue,99}, shape=rectangle, text=black, line width=1mm, minimum width=8cm, minimum height=1cm]
\tikzstyle{blockprime}=[fill={rgb,255: red,108; green,154; blue,195}, draw=black, shape=rectangle, text=white]
\tikzstyle{Hom3by3}=[fill=white, draw=black, shape=rectangle, minimum height=20mm]
\tikzstyle{finsetdot}=[fill=black, draw=black, shape=circle]
\tikzstyle{Hom1x3}=[fill=white, draw=black, shape=rectangle]
\tikzstyle{new edge style 1}=[line width=0.6mm]
\tikzstyle{Hom1x3_rounded}=[fill=white, draw=black, shape=rectangle, rounded corners]
\tikzstyle{outerbox6x17}=[fill=none, draw=black, shape=rectangle, minimum height=6cm, minimum width=17cm, rounded corners]
\tikzstyle{outerbox}=[fill=none, draw=black, shape=rectangle, rounded corners, minimum height=4cm, minimum width=6cm]
\tikzstyle{add}=[fill=white, draw=black, shape=circle]
\tikzstyle{rel_dagger}=[fill=white, draw=black, shape=trapezium, rounded corners, shape border rotate=90, minimum width=1.5cm]
\tikzstyle{rel}=[fill=white, draw=black, shape=trapezium, rounded corners, shape border rotate=270, minimum width=1.5cm]
\tikzstyle{mhom3x1}=[fill=white, draw=black, shape=rectangle, rounded corners, minimum height=10mm]
\tikzstyle{wamp}=[rounded rectangle, rounded rectangle west arc=0pt, draw, inner sep = 2pt,minimum width=5ex,minimum height=3ex,fill=white] 
\tikzstyle{wamp_dagger}=[rounded rectangle, rounded rectangle east arc=0pt, draw, inner sep = 2pt,minimum width=5ex,minimum height=3ex,fill=white] 

\tikzstyle{new edge style 0}=[->, draw={rgb,255: red,127; green,129; blue,131}, line width=0.6mm]
\tikzstyle{arrow}=[draw={rgb,255: red,177; green,179; blue,182}, ->, line width=0.4mm]
\tikzstyle{undirected}=[-, draw={rgb,255: red,177; green,179; blue,182}, line width=0.6mm]
\tikzstyle{cospan_edge}=[->, draw=black, dashed, line width=0.4mm]
\tikzstyle{undirected edge}=[-, line width=0.6mm, draw={rgb,255: red,177; green,179; blue,182}]
\tikzstyle{blueline}=[-, draw={rgb,255: red,108; green,154; blue,195}, line width=0.7mm]
\tikzstyle{functor}=[->, draw=black, line width=1.2mm]

\usepackage{assets/styles/quiver}

\usepackage{biblatex}
\usepackage{svg}

\tikzset{touch src/.style={shorten <=-4pt}}
\tikzset{touch tgt/.style={shorten >=-4pt}}
\usetikzlibrary{fit,positioning,shapes.geometric}
\usetikzlibrary{cd}
\tikzset{
  touch src/.style={shorten <=-4pt},
  touch tgt/.style={shorten >=-4pt}}
\tikzset{
  place/.style={circle,thick,draw=black!75,minimum size=5mm,inner sep=0.25mm},
  big place/.style={ellipse,thick,draw=black!75,minimum size=5mm,inner sep=0.25mm},
  transition/.style={rectangle,thick,draw=black!75},
  anon place/.style={circle,fill=black,minimum size=2.5mm,inner sep=0},
  sum place/.style={circle,thick,draw=black!75,minimum size=2.5mm,inner sep=0}}
\tikzset{ 
  touch src/.style={shorten <=-4pt},
  touch tgt/.style={shorten >=-4pt}}
\tikzset{ 
  place/.style={ellipse,thick,draw=black!75,minimum size=5mm,inner sep=0.25mm},
  big place/.style={ellipse,thick,draw=black!75,minimum size=5mm,inner sep=0.25mm},
  transition/.style={rectangle,thick,draw=black!75},
  sum place/.style={circle,thick,draw=black!75,minimum size=2.5mm,inner sep=0}}
\tikzset{ 
  boxUWD/.style={rectangle,draw,rounded corners}
}

\definecolor{gblue}{HTML}{6C9AC3}
\definecolor{gorange}{HTML}{E28F41}
\definecolor{gbrown}{HTML}{CFB691}
\definecolor{gsage}{HTML}{CFDBCB}
\definecolor{gcream}{HTML}{FFEFCF}
\definecolor{ggrey}{HTML}{B1B3B6}
\definecolor{gforest}{HTML}{5E8E3F}
\definecolor{gpurple}{HTML}{593674}
\definecolor{gyellow}{HTML}{FCAF17}
\definecolor{gyellow}{HTML}{FCAF17}
\definecolor{gred}{HTML}{D7182A}
\definecolor{gteal}{HTML}{A8DCD9}
\definecolor{ggreen}{HTML}{80BE63}
\definecolor{ggbeige}{HTML}{E0E5CD}

\definecolor{jlred}{HTML}{CB3C33}
\definecolor{jlblu}{HTML}{4063D8}
\definecolor{jlgrn}{HTML}{389826}
\definecolor{jlprp}{HTML}{9558B2}

\newcommand{\dif}{\mathrm{d}}
\newcommand{\Dif}{\mathrm{D}}
\renewcommand{\vec}[1]{\bm{#1}}
\newcommand{\Div}{\nabla\cdot}
\newcommand{\Grad}{\nabla}

\renewcommand{\vec}[1]{\bm{#1}}

\usepackage[hidelinks]{hyperref}
\newcommand*{\email}[1]{\href{mailto:#1}{\nolinkurl{#1}} } 

\newcommand{\sfield}{\mathbb{R}_\Gamma}
\newcommand{\vfield}{\mathbb{R}^2_\Gamma}

\title{Decapodes: A Diagrammatic Tool for Representing, Composing, and Computing Spatialized Partial Differential Equations}

\author[1,*]{Luke Morris}
\author[2]{Andrew Baas}
\author[3]{Jesus Arias}
\author[3]{Maia Gatlin}
\author[4]{Evan Patterson}
\author[1]{James P. Fairbanks}
\affil[1]{Department of Computer \& Information Science and Engineering, University of Florida,
Gainesville, FL, USA}
\affil[2]{Department of Mathematics, Texas State University, San Marcos, TX, USA}
\affil[3]{Georgia Tech Research Institute, Atlanta, GA, USA}
\affil[4]{Topos Institute, Berkeley, CF, USA}
\affil[*]{Corresponding Author: \email{luke.morris@ufl.edu}}

\begin{document}
\maketitle

\begin{abstract}
  We present Decapodes, a diagrammatic tool for representing, composing, and solving partial differential equations. Decapodes provides an intuitive diagrammatic representation of the relationships between variables in a system of equations, a method for composing systems of partial differential equations using an operad of wiring diagrams, and an algorithm for deriving solvers using hypergraphs and string diagrams.
  The string diagrams are in turn compiled into executable programs using the techniques of categorical data migration, graph traversal, and the discrete exterior calculus.
  The generated solvers produce numerical solutions consistent with state-of-the-art open source tools as demonstrated by benchmark comparisons with SU2. These numerical experiments demonstrate the feasibility of this approach to multiphysics simulation and identify areas requiring further development.
\end{abstract}

\section{Introduction}

In February 2013, a group of 45 researchers representing universities and national labs around the world published an article titled \textit{Multiphysics simulations: Challenges and opportunities}, assessing the current state and future of multiphysics simulation tools~\cite{keyes_multiphysics_2013}. Multiphysics systems, or systems in which multiple independent yet coupled physical phenomena are represented, have presented a significant challenge within the simulation community. The authors describe the problem faced by researchers who need to simulate multiphysics systems: the ``buy versus build'' dilemma. In order to simulate multiphysics systems, the researcher is faced with a choice between purchasing an expensive, application-specific tool or constructing their own tool, a process which could be similarly expensive and more risky. Even when purchasing a tool was an option, significant further development was often required to integrate this purchased software into the scientist's existing process.
We present in this paper a tool which addresses this dilemma and takes a step toward achieving the longer-term goals set forth~\cite{keyes_multiphysics_2013}:

\begin{quote}
    A longer-term and ambitious goal might be to allow domain experts who do not have extensive knowledge of multiphysics coupling and possibly have little advanced HPC expertise to write robust, coupled, parallel simulation codes that exhibit reasonably good performance.
\end{quote}

Our tool, called \emph{Discrete Exterior Calculus Applied to Partial and Ordinary Differential Equations} (Decapodes), enables scientists to declaratively specify physical equations and then automatically generate a simulator that solves the equations. The Decapodes software allows scientists to rapidly create new simulators as their physical modeling assumptions change. The main contribution of this work is to show how
scientific computing workflows can be improved by developing software with mathematical abstractions that connect to the mathematical models being solved. We use category-theoretic methods to formalize the physical models and simulations so that solver code can be generated from the formal specifications. 

A key problem this paper addresses is how much time it takes to develop and evolve a simulator for multiple physics. As simulators confront increasingly coupled physical systems with a greater diversity of governing equations, instabilities naturally arise. Often, these instabilities are handled using numerical methods designed for specific coupling problems, resulting in complex software which can only be understood, let alone modified, by its original developers.
Since the physical model and numerical methods are tightly intertwined, improving just one of these aspects of the simulator requires complete knowledge of both.
By automating more of this process, we can improve the productivity of scientific software developers.

Scientific software can be organized onto a spectrum of how explicitly the software represents the mathematical model it implements. On one side are computer algebra systems that have data structures representing mathematical formulas and use algorithms to manipulate these formulas to solve systems of equations. Solutions are also represented as formulas, which can be evaluated for particular choices of physical parameters. Such tools excel at providing deep understanding of the mathematical model and, insofar as the model accurately captures a physical phenomenon, they can provide deep insight into the physical phenomenon.
However, the simplifying assumptions needed to obtain closed-form solutions generally do not extend to complex multiphysics systems on complex geometries.
In order to understand such physical systems, one must turn to numerical solutions. Programs that implement solvers to compute numerical solutions lie at the other end of the spectrum. Numerical programs typically have no explicit, symbolic representation of the system and rely on knowledgeable programmers to carefully design programs that compute the correct solution for each physical situation. It is challenging to build simulators that are flexible to the choice of spatial mesh, parameter regimes, numerical method, and governing equations. Traditional scientific software development methodologies couple all of these design decisions into a single program. 

Existing frameworks such as Firedrake~\cite{10.1145/2998441} and FEniCS~\cite{logg_automated} achieve some of the longer-term goals set by Keyes et al.~\cite{keyes_multiphysics_2013} by using the Unified Form Language (UFL)~\cite{10.1145/2566630} and generating solvers from UFL expression trees. As Rathgeber et al.~\cite{10.1145/2998441} note, FEniCS introduced a ``separation of concerns'' betweeen ``$\textit{employing}$ the finite element method and $\textit{implementing}$ it,'' while Firedrake introduces a ``separation within the implementation layer.'' However, the compositions of equations formulated via UFL are not themselves treated as first-class objects. So, reusing previously-written equations still requires the manual editing of those equations, which leaves a large surface area for bugs to be introduced.

In simulators that lack a symbolic representation of the physical system such as SU2~\cite{EconomonSU2}, small changes to the governing equations cause cascading changes to the software requirements. The choice of discretization scheme, numerical methods, operator implementations, and meshing tools can depend on the choice of governing equations. This might not be a problem if the full formulations of physical laws were always used, but in building high-performance simulators, scientists need to use engineering approximations that allow for larger systems to be solved in finer detail than would otherwise be feasible. While the fundamental laws of nature might not change often or at all, the approximations used by computational scientists implementing simulators do change rapidly. Or at least they would change rapidly if the technology used to implement the simulators made such changes quick and easy.

\begin{figure}
    \centering
    \begin{subfigure}{0.3\textwidth}
    \includegraphics[width=\textwidth]{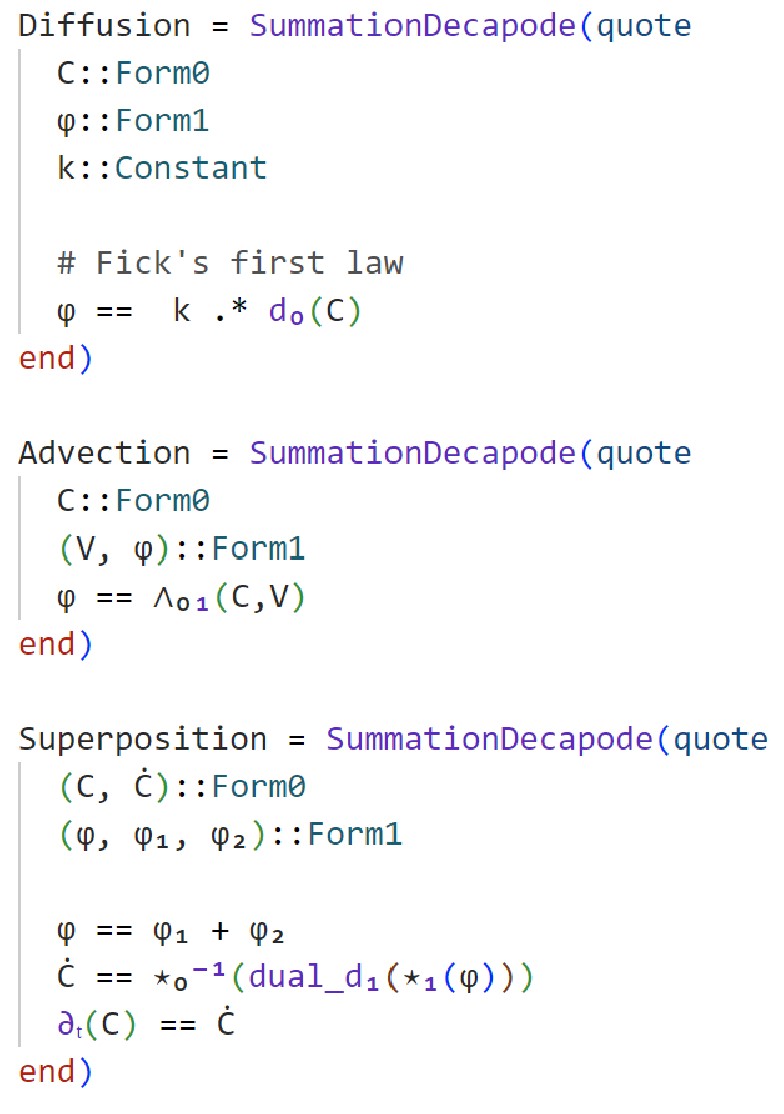}
    \label{fig:julia-blocks}
    \caption{Discrete Exterior Calculus Equations.}
    \end{subfigure}
    \begin{subfigure}{0.4\textwidth}
    \begin{subfigure}{0.25\textwidth}
    \includegraphics[width=\textwidth]{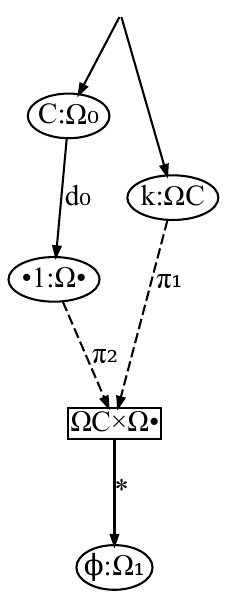}
    \label{fig:snapshot-diffusion}
    \end{subfigure}
    \begin{subfigure}{0.25\textwidth}
    \includegraphics[width=\textwidth]{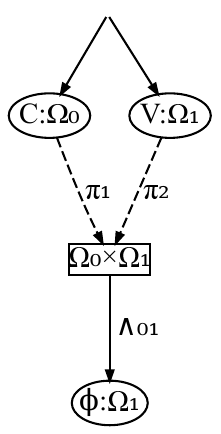}
    \label{fig:advection}
    \end{subfigure}
    \begin{subfigure}{0.33\textwidth}
    \includegraphics[width=\textwidth]{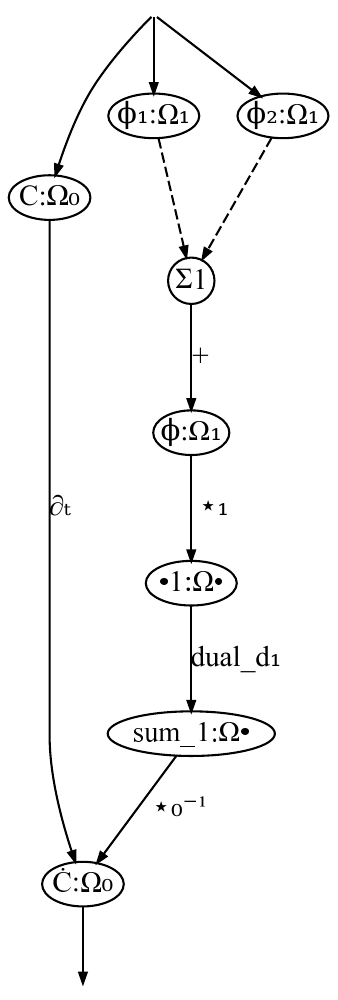}
    \label{fig:superposition}
    \end{subfigure}
    \caption{Formal Multiphysics Diagrams.}
    \end{subfigure}
    
    \begin{subfigure}{0.3\textwidth}
    \includegraphics[width=\textwidth]{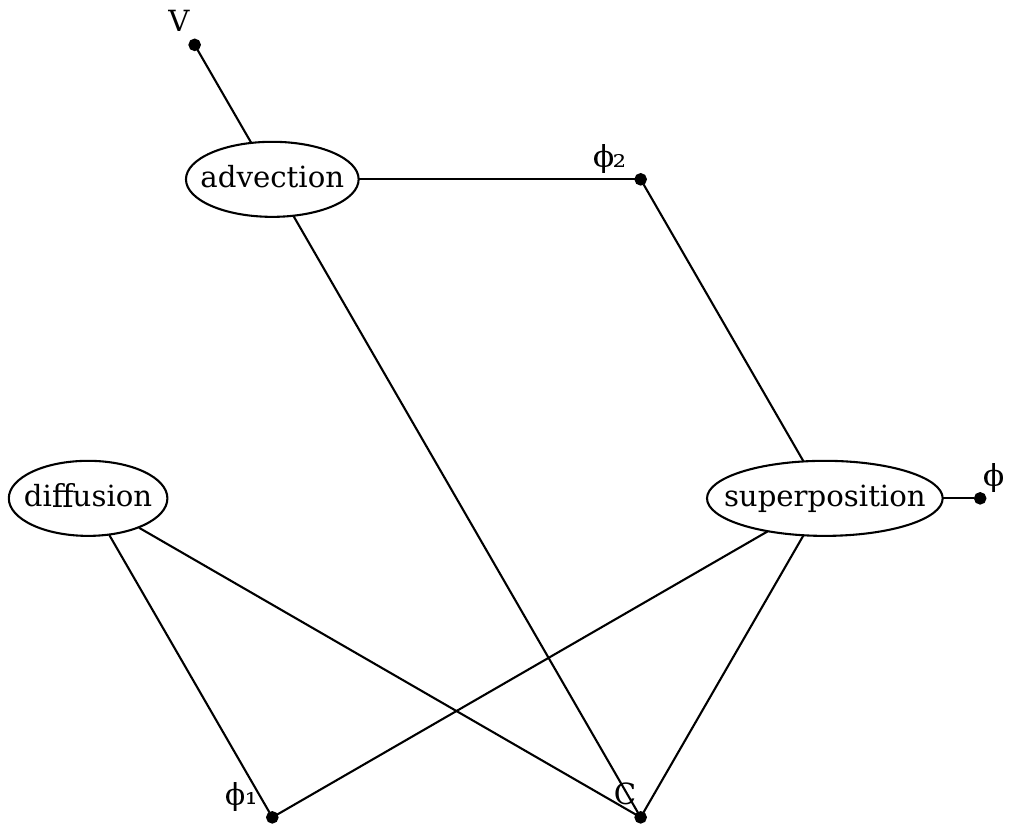}
    \caption{Operadic Compositions.}
    \label{fig:compose-diffadv}
    \end{subfigure}
    \begin{subfigure}{0.3\textwidth}
    \includegraphics[width=\textwidth]{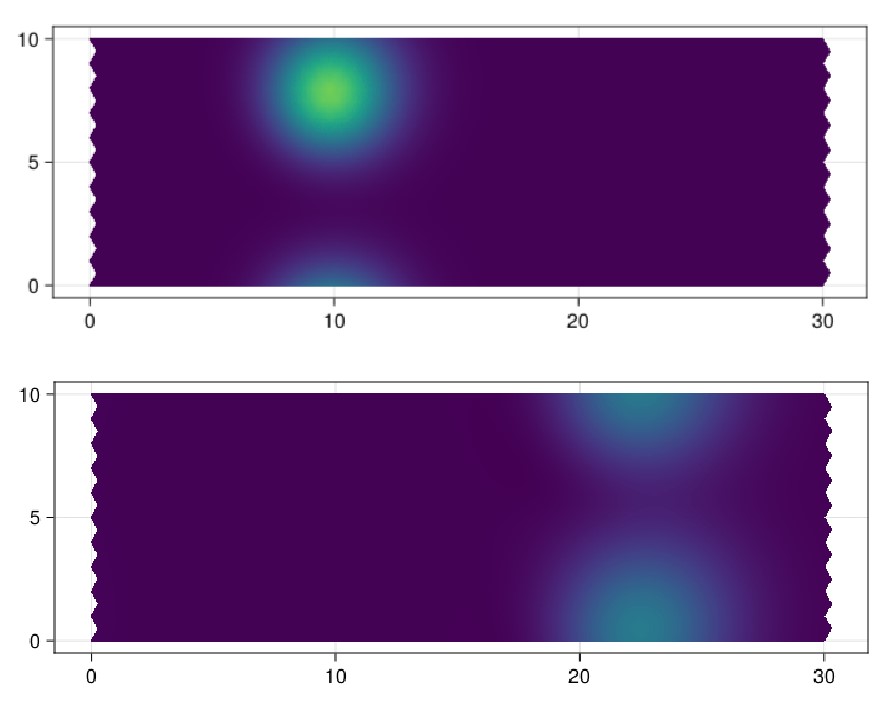}
    \caption{Executable Simulators.}
    \label{fig:diffadv-sim-results}
    \end{subfigure}
    \caption{A snapshot of the Decapodes workflow.}
\end{figure}

We address these challenges by introducing two graphical formalisms to the field of multiphysics simulation: diagrams of elements to describe governing equations, and directed hypergraphs to describe solver computations. Both formalisms are modular, compositional, and generalizable over many domains of mathematics. Using these graphical languages, we propose and implement a computational system for multiphysics simulation which captures the structure of both the mathematical description of systems and the resulting computation for those systems. Our system is structured enough to provide an intuitive interface for the end-user and flexible enough to allow for the introduction of optimizations at several levels.

The thesis of this paper is that using the Decapodes tooling of multiphysics diagrams provides the following benefits to the computational physicist. Encoding a multiphysics system with Decapodes provides an automatically-generated simulation that is accurate and generalizes over choice of mesh. One can construct more complicated multiphysics systems from simpler Decapodes using composition, thus using composition patterns to manage the complexity of large multiphysics systems. By applying the paradigm of hierarchical composition in the field of computational physics, we take a step toward achieving the \textit{Multiphysics simulations} vision of a truly collaborative and unified scientific modeling community.

\subsection{Tonti Diagrams and Diagrammatic Equations}

Diagrams are routinely used in science to provide an intuitive representation of a physical system or mathematical model. Feynman diagrams, circuit diagrams, free body diagrams, and bond graphs, among many others, are graphical tools that aid the understanding of complex and coupled phenomena.

The utility of these diagrammatic representations depends on their degree of precision. Feynman diagrams and circuits have established conventions by which one can read a diagram and understand the described system, thus providing a common language for scientists working in these areas. If a diagrammatic representation is sufficiently formal, an algorithm can be implemented that takes the diagram as input and emits a mathematical model expressed in traditional equational notation. Issues arise when informal conventions are introduced into these diagrams, such as a scientist adding descriptive notes to a circuit diagram that invalidate the standard interpretation or an instructor adding their own conventions to free body diagrams. Despite having extra flexibility, an informal graphical representation becomes confusing and breaks down as a reliable medium for communication. A modern approach to rigorous diagrams is exemplified by Censi's Co-Design Diagrams, which provide a rigorous language for describing the design of complex systems as optimization problems~\cite{censi2016codesign}. Fong and Spivak present a collection of formal diagrammatic languages for mathematics in teaching the techniques of applied category theory~\cite{Fong_Spivak_2018}. Our approach builds on the categorical foundations of databases~\cite{spivak2014} to represent mathematical objects as instances of relational databases~\cite{patterson_categorical_2022}.

We are not the first to use diagrams to describe the system of equations defining a physical theory. Diagrams involving physical quantities and their relations appear throughout the
literature on computational physics and the exterior calculus, sometimes called \emph{Tonti
diagrams} {\cite{tonti1972,tonti_mathematical_2013}} or, in the context of electromagnetism,
\emph{Maxwell's house} {\cite{deschamps1981,bossavit1998d}}.
Our approach is inspired by the work of Enzo Tonti, a physicist who was a major proponent of diagrams as tools for describing physical systems. Tonti developed a diagrammatic language based on the structure of the exterior calculus and encoded both variables and their relationships in a physical system~\cite{tonti_mathematical_2013}. An example of a Tonti diagram is shown in Fig.~\ref{fig:tonti_ns}. In this diagram, variables are the vertices of a directed graph and operators acting on the variables are directed edges.

\begin{figure}
\centering
\begin{tiny}
\[\begin{tikzcd}
	\arrow["{L^h_k = \nabla_k v^h}"{description}, from=4-1, to=8-1]
	\arrow[dash, from=8-1, to=7-4]
	\arrow[dash, from=first_split.center, to=6-6]
	\arrow[dash, from=7-4, to=7-6]
	\arrow[dashed, dash, from=6-6, to=lambda_eq]
	\arrow[dashed, dash, from=7-6, to=mu_eq]
	\arrow[dashed, dash, from=lambda_eq, to=first_join.center]
	\arrow[dashed, dash, from=mu_eq, to=7-10]
	\arrow[dash, from=6-10, to=second_join.center]
	\arrow[from=7-10, to=7-12]
	\arrow["{p_i \stackrel{\text{mat}}{=} g_{ih}\rho v^h}"{description}, from=4-1, to=4-13]
	\arrow[dash, from=8-1, to=8-13]
	\arrow[from=8-13, to=4-13]
	\arrow[from=5-2, to=8-1]
	\arrow[from=1-2, to=4-1]
	\arrow[from=8-13, to=5-14]
	\arrow[from=7-12, to=5-14]
	\arrow[dash, from=5-2, to=5-14]
	\arrow[dash, from=1-2, to=1-14]
	\arrow[from=5-14, to=1-14]
	\arrow["{\partial_t p_i + \nabla_h (p_i v^h) - \nabla_h \tau^h_i = f_i}"{description}, from=4-13, to=1-14]
	\arrow[from=1-2, to=5-2]
 & \square &&&&&&&&&&&& {f_i} \\
	\\
 &&&&&&&&&&&&&& \\
	{v^h} &&&&&&&&&&&& {p_i} \\
	& \square &&&&&&&&&&&& {t^h_i} \\
	&&&&& \theta && \node (lambda_eq) [fill=lightgray] {(\lambda+\frac{2}{3}u)\theta\delta_{ik}}; && {-p \delta_{ik}} \\
	&&& {D_{ik}} & \node (first_split) {}; & {D^\prime_{ik}} && \node (mu_eq) [fill=lightgray] {2\mu D^\prime_{ik}}; & \node (first_join) {}; & {\sigma^\prime_{ik}} & \node (second_join) {}; & {\sigma_{ik}} \\
	{L^h_k} &&&&&&&&&&&& \square
\end{tikzcd}\]
\end{tiny}
    \caption{This recreation of Diagram FLU6-12 \cite{tonti_mathematical_2013} from Tonti's 2013 work provides a graphical representation of the Navier-Stokes equations. Nodes labeled as $\square$ do not have relevance to the physics, but are necessary to represent the geometric layout of nodes, which encodes whether quantities are in primal or dual space and time. One equation is shared between two distinct edges, and some are color-coded with dashed arrows.}
     \label{fig:tonti_ns}
\end{figure}

Benefits of the Tonti diagram formalism include:

\begin{enumerate}
    \item Ease of understanding: relationships between quantities in a physical model are more legible, potentially allowing non-expert involvement in the interrogation of the model.
    \item Use of exterior calculus: Tonti diagrams integrate naturally with the exterior calculus (as described by Tonti in his 2013 book~\cite{tonti_mathematical_2013}) and so naturally suggest discretization schemes and numerical simulation processes.
    \item Connection with solution methods: Tonti diagrams show the interdependence of variables in the model which an expert can use to construct an iterative solver for the system.
    \item Compositionality: the composition of two physical models can be intuitively understood as identifying shared variables.
\end{enumerate} 

While Tonti diagrams provide a useful interface for describing physical theories, they begin to break down into informality as the physics they describe becomes more complex.
For example, Tonti's description of the Navier-Stokes equations (Fig.~\ref{fig:tonti_ns}) introduces new graphical notations and typed notes in order to communicate the equations.
While bespoke visual and textual cues may suffice for a scientist to extract the intended meaning behind a particular diagram, they demonstrate a clear need to extend and clarify Tonti diagrams if they are to be used as data structures in a computational system. Unfortunately, the framework of Tonti diagrams is insufficiently rigorous to provide a foundation for multiphysics systems due to a number of problems:

\begin{enumerate}
    \item Saturation: the diagrams can saturate. Since the location of a variable is determined by its dimensionality and primality, there can exist only one primal 0-form, primal 1-form, etc. Multiphysics diagrams contain too many variables to represent with fixed vertex locations.
    \item Ambiguity: multiple arrows incoming to a node sometimes instruct the reader to apply a pair of constraints simultaneously and sometimes to add forms. The reader must use their physical knowledge to resolve this ambiguity.
    \item Restriction to unary operations: multivariate operations such as wedge products and Lie derivatives cannot be shown on the Tonti diagram, which limits their applicability to fluid mechanics models.
\end{enumerate}

These problems with Tonti diagrams are resolved by Patterson et al. through the construction of categories of diagrams of elements~\cite{patterson2022diagrammatic}.
These categories offer both a syntax for combinatorial diagrams (directed labeled graphs) and a semantics for systems of equations.
Diagrams are viewed as presenting small categories, and functors from these categories to analytic categories---such as those of vector spaces, or of cochains over a manifold---are semantics that allow the diagrams to represent equations.
The authors provide a theoretical framework for understanding how the combinatorial properties of the diagrams can be used to study the properties of the solutions of equations.

Though formalized by applied category theory, the graphical nature of the tools allows interpretation by practicing mathematicians, scientists and engineers without training the categorical methods used to develop the tools. Any scientist can understand the diagrammatic equations presented below. Variables are drawn as nodes labeled $x:X$, meaning that $x$ is an element of $X$, and equations are drawn as arrows $x:X\xrightarrow{f} y:Y$, meaning that $y$ is the result of applying the operator $f: X \to Y$ to $x$. Circles containing $\Sigma$ are used to encode $n$-ary summation, and boxes are used to to encode multilinear operators.

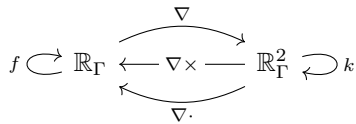
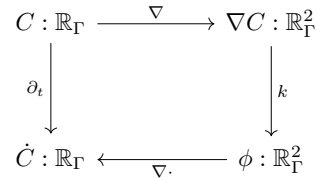
\begin{figure}[htbp]
    \centering
    \begin{subfigure}[b]{0.4\textwidth}
        \centering
\[\begin{tikzcd}
	{\mathbb{R}_\Gamma} && {\mathbb{R}^2_\Gamma}
	\arrow["\nabla", shift left=2, curve={height=-8pt}, from=1-1, to=1-3]
	\arrow["\nabla\cdot", shift left=2, curve={height=-8pt}, from=1-3, to=1-1]
	\arrow["\nabla\times"{description}, from=1-3, to=1-1]
    \arrow["f", loop left, from=1-1, to=1-1]
    \arrow["k", loop right, from=1-3, to=1-3]
\end{tikzcd}\]
        \caption{Theory-level diffusion diagram}
        \label{fig:diff_theory}
    \end{subfigure}
    \hspace{0.1\textwidth}
    \begin{subfigure}[b]{0.4\textwidth}
        \centering
        \adjustbox{scale=0.93,center}{
\begin{tikzcd}
	{C:\mathbb{R}_\Gamma} && {\nabla C:\mathbb{R}^2_\Gamma} \\
	\\
	{\dot C: \mathbb{R}_\Gamma} && {\phi:\mathbb{R}^2_\Gamma}
	\arrow["\nabla", from=1-1, to=1-3]
	\arrow["\nabla\cdot", from=3-3, to=3-1]
	\arrow["k", from=1-3, to=3-3]
	\arrow["{\partial_t}"', from=1-1, to=3-1]
\end{tikzcd}
        }
        \caption{Diffusion equations expressed as a diagram}
        \label{fig:diff_diagram}
    \end{subfigure}
    \caption{The theory of two-dimensional vector calculus includes the traditional operators of multivariate calculus (\ref{fig:diff_theory}). $\R_\Gamma$ and $\R^2_\Gamma$ are the spaces of scalar and vector fields over a domain $\Gamma \subseteq \R^2$. The operators $\nabla, \nabla\cdot, \nabla\times$ are the gradient, divergence, and curl, while $f$ and $k$ represent pointwise functions applied everywhere over $\Gamma$. These operators provide a language in which to express models such as the diffusion equation (\ref{fig:diff_diagram}). Variables and operators in the model are typed with spaces and maps in the theory. 
    }
    \label{fig:diffusion}
\end{figure}

As can be seen from Fig.~\ref{fig:diffusion}, these diagrams are similar to Tonti diagrams. The variables involved in the diffusion equation are presented as vertices in Fig.~\ref{fig:diff_diagram}, and the operators to calculate one variable from another are present on the edges. Unlike in Tonti diagrams, if two edges point to the same vertex, this unambiguously imposes an equality between the results of the operations of those edges. The diagram in Fig.~\ref{fig:diffusion} states that
\begin{equation}\label{eq:c_dot}
    \dot{C} = \partial_t C 
\end{equation}
and
\begin{equation}
    \dot{C} = \Div k \Grad C,
\end{equation}
where each path in the diagram represents an equation and equations are coupled by shared variables. Note that $\dot{C}$ is here a purely formal symbol that attains its usual meaning from eq. \ref{eq:c_dot}.

A key benefit of the categorical approach to diagrammatic equations is the ability to construct higher level operations on these diagrams. For example, Patterson et al provide a language for specifying multiphysics models by composing diagrams by gluing vertices together. This composition builds on the work of Baez et al.~\cite{Baez_2022}~\cite{Baez_Pollard_2017}~\cite{Baez_Fong_2018} by describing multiphysics models as structured cospans of diagrams. In essence, models compose by sharing variables.

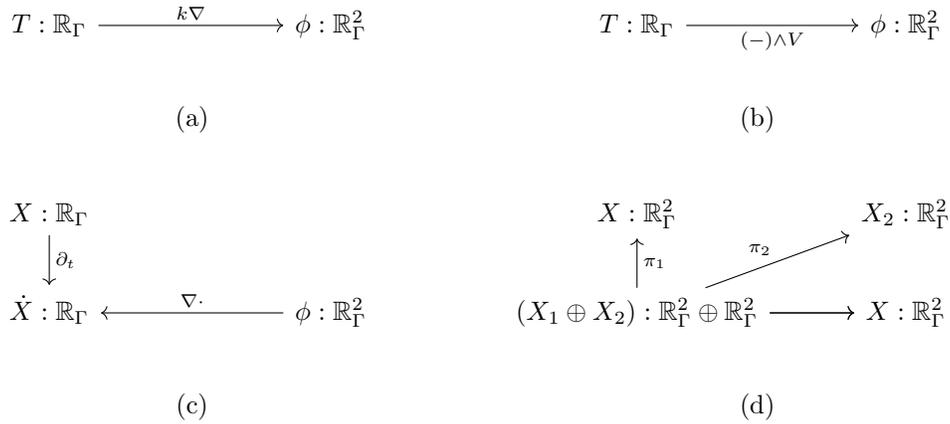
\begin{figure}[htbp]
    \centering

\[\begin{tikzcd}
	{T:\R_\Gamma} && {\phi:\R^2_\Gamma} && {T:\mathbb{R}_\Gamma} & {\phi:\R^2_\Gamma} \\
	& {\text{(a)}} &&& \node [xshift=16mm] {\text{(b)}}; \\
	{X:\R_\Gamma} &&&& {X:\R^2_\Gamma} & {X_2:\R^2_\Gamma} \\
	{\dot X: \R_\Gamma} && {\phi:\R_\Gamma^2} && {(X_1\oplus X_2):\mathbb{R}^2_\Gamma\oplus\R_\Gamma^2} & {X:\R^2_\Gamma} \\
	& {\text{(c)}} &&& \node [xshift=16mm]  {\text{(d)}};
	\arrow["\nabla\cdot"', from=4-3, to=4-1]
	\arrow["{\partial_t}", from=3-1, to=4-1]
	\arrow["k\nabla", from=1-1, to=1-3]
	\arrow["{(-)\wedge V}"', from=1-5, to=1-6]
	\arrow["{\pi_1}"', from=4-5, to=3-5]
	\arrow["{\pi_2}", from=4-5, to=3-6]
	\arrow[from=4-5, to=4-6]
	\arrow[from=4-5, to=4-6]
\end{tikzcd}\]
    \caption{The four physical components used to define the advection diffusion system. The components are (a) Fick's law, (b) advection of a scalar field along a vector field, (c) conservation of mass, and (d) linear superposition of flux. }
    \label{fig:components}
\end{figure}

\begin{figure}[htbp]
\begin{subfigure}[b]{0.9\textwidth}
\centering
\begin{tikzpicture}[baseline=(current bounding box.center),tips=proper]
        \node (C) {$T: \sfield$};
        \node (F) [below=4em of C] {$\phi: \vfield$};
        \node (CF) [fit=(C) (F), inner sep=0] {};
        \node (superposition) [boxUWD, left=of F, align=center] {Flux \\ superposition};
        \node (F1) [above=4em of superposition] {$\phi_1: \vfield$};
        \node (F2) [left=2em of superposition] {$\phi_2: \vfield$};
        \node (advection) [boxUWD, above=2em of F2] {Advection};
        \node (diffusion) [boxUWD, above=of advection] {Diffusion};
        \node (conservation) [boxUWD, right=of CF, align=center] {Mass \\ conservation};
        \node (U) [left=1em of advection] {$u: \vfield$};
        \node (outer) [boxUWD, inner sep=1em, fit=(C) (F) (diffusion) (advection) (conservation) (U)] {};
        \path
            (diffusion.east) edge[out=0, looseness=0.5] (C.north west)
            (diffusion.south) edge[out=270, in=180] (F1.west)
            (advection.south) edge (F2)
            (conservation.165) edge[out=180, in=0] (C)
            (conservation.195) edge[out=180, in=0] (F)
            (superposition) edge (F)
            (superposition) edge (F1)
            (superposition) edge (F2)
            (outer.north) edge[out=270] (C.north)
            (outer.west) edge[in=180, out=0] (U.west)
            (U.east) edge[in=180, out=0] (advection.west);
        \path[draw,line width=5pt,white]
            (advection.east) edge[out=0, in=180] (C.west);
        \path
            (advection.east) edge[out=0, in=180] (C.west);
    \end{tikzpicture}

\caption{}\label{fig:pattern}
\end{subfigure}
\centering
    \begin{subfigure}[b]{0.90\textwidth}
    \centering
\[\begin{tikzcd}
	{T:\R_\Gamma} && {\phi_2:\R^2_\Gamma} \\
	\\
	& {\phi_1:\R^2_\Gamma} & {(X_1\oplus X_2):\mathbb{R}^2_\Gamma\oplus\R_\Gamma^2} \\
	{\dot T: \R_\Gamma} && {\phi:\R_\Gamma^2}
	\arrow["\nabla\cdot"', from=4-3, to=4-1]
	\arrow["k\nabla", from=1-1, to=3-2]
	\arrow["{\partial_t}"', from=1-1, to=4-1]
	\arrow["{\pi_2}", from=3-3, to=1-3]
	\arrow["{\pi_1}"', from=3-3, to=3-2]
	\arrow["{+}", from=3-3, to=4-3]
	\arrow["{(-)\wedge V}", from=1-1, to=1-3]
\end{tikzcd}\]
        \caption{}
        \label{fig:adv_diff_composition}
    \end{subfigure}
    \caption{An undirected wiring diagram describing how to connect physical principles of diffusion, advection, flux superposition, and mass conservation to create a composite advection-diffusion system (\ref{fig:pattern}) as given by Patterson et al. \cite{patterson2022diagrammatic}. This pattern can be applied to the component physics systems in Fig~\ref{fig:components} to create the composite advection-diffusion system (\ref{fig:adv_diff_composition}). Note that the operation relating $T$ and $\phi_2$ is the wedge product parameterized by some fixed vector field $V$, similar to a parameterized Lie derivative.}\label{fig:composition}
\end{figure}
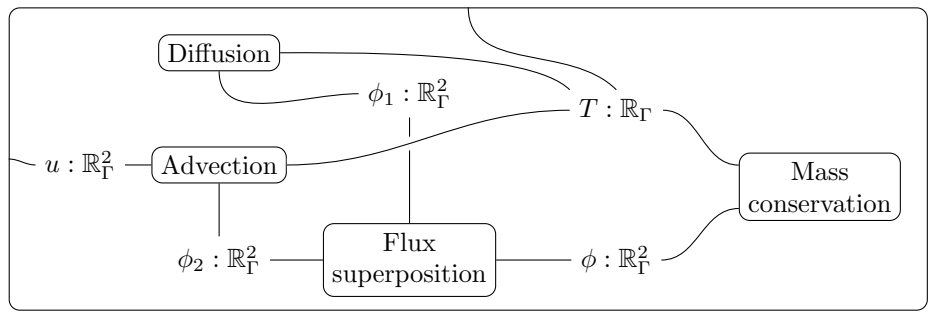
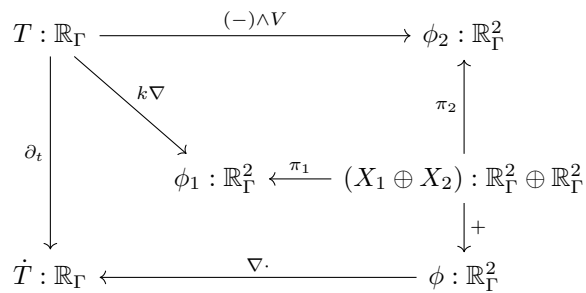

Another benefit of using diagrams of elements for describing physical systems is that the diagrams can be composed hierarchically using the operad of undirected wiring diagrams (UWDs).
A UWD consists of boxes, ports, junctions, and wires. Ports are attached to boxes, and wires connect ports to junctions.
UWDs can act as patterns of composition for individual diagrammatic equations. An example of this composition is provided in Fig.~\ref{fig:adv_diff_composition}. Taking the individual physics described by the diagrams in Fig.~\ref{fig:components} and composing them with the UWD in Fig.~\ref{fig:pattern}, we are able to generate the composed physical system with advection, diffusion, and mass conservation in Fig.~\ref{fig:adv_diff_composition}. In this way, simpler component physics can be composed into significantly more complex multiphysics. A treatment of the mathematical and theoretical aspects of this work can be found in Patterson et al.~\cite{patterson2022diagrammatic}.

Having specified a PDE with a diagrammatic equation, one can construct an ODE to solve it via the method of lines, which in turn we solve with an explicit time stepping method such as Runge-Kutta. We seek a systematic approach to representing both the diagrammatic equation and the corresponding vector field as combinatorial objects. For the vector field we turn to a traditional approach of representing computations as graphs.

\subsection{Computation Graphs and the Method of Lines}

We represent a physical model by a system of partial differential equations and solve time-varying problems with the method of lines to reduce it to a system of ordinary differential equations (ODEs). The right hand side of that ODE is given by a hypergraph representing the computation of $\dot{u} = f(u)$.

While Tonti diagrams and diagrammatic equations encode the constraints that a physical system must satisfy, computation graphs encode processes for computing updates to the state of the physical system. In a typical approach, one represents each variable in the system with a node in a dependency graph and adds an edge from vertex $x$ to vertex $y$ if the computation of $y$ depends on the value of $x$. Dependency graphs help programmers understand how information flows through a computation. However, dependency graphs under-specify the computation because there is no representation of the functions in the graph. In this sense, dependency graphs are like Tonti diagrams because they rely on an intelligent interpreter to derive a specific execution of the dependency graph.

An alternative to dependency graphs is to use directed hypergraphs to represent a computation.  Each variable in the system is a vertex in the hypergraph and each hyperedge has an ordered tuple of input vertices and an ordered tuple of output vertices. This encoding of a computation uses directed hyperedges to represent the multiple-input-multiple-output functions common in mathematical computing. For example the assignment $(y_1, y_2) = f(x_1,x_2,x_3)$ is represented as a hyperedge with two outputs and three inputs.
Hypergraphs ``give a canonical, combinatorial representation for string diagrams'' as arbitrary symmetric monoidal categories~\cite{bonchi2022string} (Sec 3). This formalizes and generalizes the graphical depiction of computation graphs used to describe data-flow programs~\cite{conal2017compiling}.

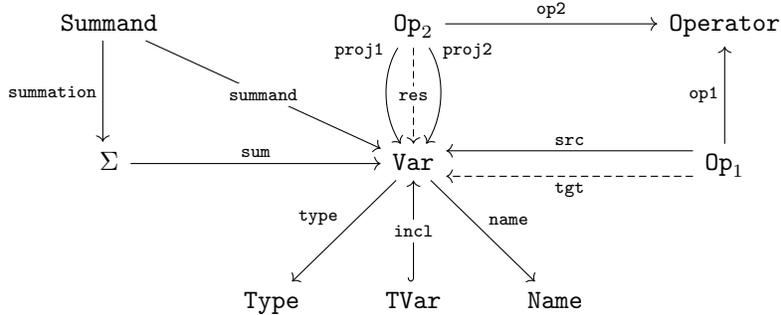
\begin{figure}
\[\begin{tikzcd}
	\texttt{Summand} && {\texttt{Op}_2} && \texttt{Operator} \\
	\\
	\Sigma && \texttt{Var} && {\texttt{Op}_1} \\
	\\
	& \texttt{Type} & \texttt{TVar} & \texttt{Name}
	\arrow["\texttt{proj1}"'{pos=0.2}, curve={height=12pt}, from=1-3, to=3-3]
	\arrow["\texttt{res}"{description}, dashed, from=1-3, to=3-3]
	\arrow["\texttt{proj2}"{pos=0.2}, curve={height=-12pt}, from=1-3, to=3-3]
	\arrow["\texttt{src}"', shift right=2, from=3-5, to=3-3]
	\arrow["\texttt{tgt}", shift left=2, dashed, from=3-5, to=3-3]
	\arrow["\texttt{op1}", from=3-5, to=1-5]
	\arrow["\texttt{op2}", from=1-3, to=1-5]
	\arrow["\texttt{sum}", from=3-1, to=3-3]
	\arrow["\texttt{summand}"{description}, from=1-1, to=3-3]
	\arrow["\texttt{summation}"', shift right=1, from=1-1, to=3-1]
	\arrow["\texttt{type}"', from=3-3, to=5-2]
	\arrow["\texttt{name}", from=3-3, to=5-4]
	\arrow["\texttt{incl}"{description}, hook, from=5-3, to=3-3]
\end{tikzcd}\]
\caption{The Decapode schema.}
\label{fig:decapode_schema_fig}
\end{figure}

If the computation contains only unary operators, binary operators, and sums, which is common in descriptions of physical models, then the computations can be stored in a database with the schema given in Fig. \ref{fig:decapode_schema_fig}.
In a database with this schema, the table $\texttt{Var}$ stores the vertex set as rows, and the tables $\texttt{Op}_1$, $\texttt{Op}_2$ store the occurrences of unary and binary operations, respectively. The vertices in $\texttt{Var}$ are the variables of the system. The inputs and outputs of each operator are a foreign key constraint that ensures that the row corresponding to an operator points to the variables that it acts on. For example a binary operator such as the wedge product, $\wedge$, will have two input variables $i_1, i_2$ and one output variable $o_1$. Attributes of types $\texttt{Type}$ and $\texttt{Name}$ store the type and name of a variable. Rows in the $\Sigma$ table are summations, which has a variable as its sum and two or more variables as its summands. The $\texttt{TVar}$ table records which variables are tangent variables, which are the difference between the value of a variable at one time step and the next. Bonchi et al.~\cite{bonchi2022string} provides a more general construction for arbitrary monoidal signatures, but we will focus on signatures containing only binary and unary operators.

\subsection{Discrete Exterior Calculus}
\label{sec:dec}

While the specification of models and generation of vector fields is independent of any discretization scheme, the choice of discretization scheme is essential to generating simulators with good performance. We turn to the discrete exterior calculus because it is a natively discrete version of calculus for modeling spatial phenomena.

In his PhD thesis, Anil Hirani demonstrates that a discrete form of exterior calculus can be constructed entirely from simple geometric operations on simplicial complexes~\cite{hirani2003}. As suggested by the name, our Decapodes tool uses an implementation of the DEC as its numerical backend. In this subsection we summarize the features of the DEC that are important for understanding the simulation examples provided in Section~\ref{sec:experiments}. For a more complete description of the DEC, see Hirani's original thesis~\cite{hirani2003} or Crane's lecture notes on discrete differential geometry~\cite{crane_discrete_2022}.

The DEC, as a variant of the exterior calculus for discrete spaces, provides an interface for typing physical quantities by their dimensions and orientations. This scheme is closely related to Tonti's classification of physical quantities~\cite{tonti_mathematical_2013} and so provides an appropriate type system for our application.

The DEC distinguishes between values associated with different dimensionalities of mesh components, these values being called \textit{forms}. In a 2D mesh, 0-forms (zero-dimensional forms), 1-forms and 2-forms are associated with values on points, lines, and surfaces, respectively. For example, a scalar field representing concentration would be discretized as a 0-form over a mesh, where each point on the mesh has a corresponding concentration value. This results in a vector of values where each element corresponds to a spatial element of the mesh. In addition to these \textit{primal} forms, the DEC includes for each primal form a corresponding dual form, defined over a dual mesh. The dual mesh is constructed as a mesh which has $k$-dimensional elements which correspond directly to the $n-k$ dimensional elements of the original $n$-dimensional mesh, referred to as the primal mesh. In two dimensions, dual points, lines and surfaces correspond to primal surfaces, lines and points. There are many valid dual meshes for a single primal mesh, and so the process of generating a dual mesh for a primal mesh involves the choice of how dual points will be located with respect to their corresponding primal surfaces. The circumcentric and barycentric methods are two choices which we have implemented and their respective benefits and challenges are discussed in \ref{app:dec_mod}.

\begin{figure}[hbtp]
    \centering
    \begin{tikzcd}
	{\Omega_0} && {\Omega_1} && {\Omega_2} \\
	\\
	{\tilde\Omega_{2}} && {\tilde\Omega_{1}} && {\tilde\Omega_{0}}
	\arrow["d_0", from=1-1, to=1-3]
	\arrow["d_1", from=1-3, to=1-5]\
	\arrow["\star_0"', shift right=2, from=1-1, to=3-1]
	\arrow["\star_1"', shift right=2, from=1-3, to=3-3]
	\arrow["\star_2"', shift right=2, from=1-5, to=3-5]
	\arrow["{\tilde d_0}", from=3-5, to=3-3]
	\arrow["{\tilde d_1}", from=3-3, to=3-1]
	\arrow["{\star_0^{-1}}"', shift right=2, from=3-1, to=1-1]
	\arrow["{\star_1^{-1}}"', shift right=1, from=3-3, to=1-3]
	\arrow["{\star_2^{-1}}"', shift right=1, from=3-5, to=1-5]
\end{tikzcd}
    \caption{The 2D de Rham Complex}
    \label{fig:dRC}
\end{figure}
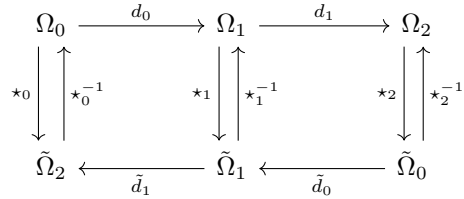

DEC operators are parameterized by the dimensionality of the form they act on, and they transform forms between different types and dimensionalities. The discrete de Rham complex shown in Figure~\ref{fig:dRC} provides a graphical depiction of the different forms and fundamental operators acting on them. Primal $n$-forms belong to the space denoted $\Omega_n$ and dual $n$-forms to the space $\tilde{\Omega}_n$. Each operator is presented as an arrow, where the domain and codomain are the source and target of the arrow, respectively. The de Rham complex shows that discrete derivatives ($d_n$ and $\tilde{d}_n$) increase the dimensionality of the forms they operate on by one, and that Hodge star operators ($\star_n$ and $\star_n^{-1}$) transform between primal and dual forms.
By selectively composing these discrete derivative and Hodge star operators, we can reproduce the traditional gradient, divergence, and curl operators from vector calculus in an entirely discrete setting. In two dimensions, the formulas are:
\begin{align}
    \nabla &= \star_1\circ d_0\\
    \nabla\cdot &= \star^{-1}_0\circ\tilde{d_1}\\
    \nabla\times &= \star_2\circ d_1.
\end{align}
This translation should enable a reader familiar with vector calculus to reason about the behavior of DEC models in the rest of this paper. For example the Laplacian for 0-forms (discrete scalar fields) is given by the familiar identity $L = \nabla\cdot \nabla = \star^{-1}_0\circ \tilde{d_1} \circ \star_1\circ d_0$.

Fig~\ref{fig:diff_diagram_dec} shows the analogous diagram for the diffusion equation in Fig~\ref{fig:diff_diagram} after making the substitutions for the discrete exterior calculus operators. 
\begin{figure}[htbp]
    \centering
    \begin{subfigure}[b]{0.4\textwidth}
        \centering
\begin{tikzcd}
  {\Omega_0} & {\Omega_0} & {\tilde\Omega_0} \\
  {\Omega_1} && {\tilde\Omega_1}
  \arrow["{\partial_t}", from=1-1, to=1-2]
  \arrow["d"', from=1-1, to=2-1]
  \arrow["k\star"', from=2-1, to=2-3]
  \arrow["{\tilde d}"', from=2-3, to=1-3]
  \arrow["{\star^{-1}}"', from=1-3, to=1-2]
\end{tikzcd}
        \caption{Theory-level diffusion diagram}
        \label{fig:diff_theory_dec}
    \end{subfigure}
    \hspace{0.1\textwidth}
    \begin{subfigure}[b]{0.4\textwidth}
        \centering
        \adjustbox{scale=0.93,center}{
\begin{tikzcd}
  C:\Omega_0 & {\dot C:\Omega_0} & \bullet:\tilde{\Omega}_0 \\
  \bullet:\Omega_1 & {} & \phi:\tilde{\Omega}_1
  \arrow["d"', from=1-1, to=2-1]
  \arrow["{\partial_t}", from=1-1, to=1-2]
  \arrow["k\star"', from=2-1, to=2-3]
  \arrow["{\tilde d}"', from=2-3, to=1-3]
  \arrow["{\star^{-1}}"', from=1-3, to=1-2]
\end{tikzcd}
        }
        \caption{Diffusion equations expressed as a diagram}
        \label{fig:diff_diagram_dec}
    \end{subfigure}
    \caption{These diagrams encode the same physical models as Fig. \ref{fig:diff_theory} and Fig. \ref{fig:diff_diagram} but in the formalism of exterior calculus instead of vector calculus. Fig.~\ref{fig:diff_theory_dec} shows the projection onto the de Rham complex.}
    \label{fig:diffusion_dec}
\end{figure}

\section{Theoretical contributions}
\label{sec:main}

In order to implement a simulator from a diagrammatic equation, we apply the method of lines to the PDE to get a vector field. The computation of this vector field is represented as a hypergraph which is then passed to an ODE solver as the right-hand side. 
This compilation process requires the Decapode to satisfy a specific structural requirement: all non-state variables must appear alone on the side of at least one equation. This equation acts as the definition of and method for computing that value. This corresponds to the property that every variable in a program must be defined to be used.

\subsection{Hypergraph Approach}\label{sec:hypergraph_alg}

In order to take advantage of category theory's existing formalism for equations-as-diagrams, we need to describe physical values as \textit{generalized elements} of some category $\cat{C}$. This concept of generalized elements is very similar to the concept of $\cat{C}$-sets presented by Patterson et al.~\cite{patterson_categorical_2022}, which essentially ``types'' each physical variable (objects of $\El(\cat{C})$) with some object in $\cat{C}$.

Diagrams in the category $\El(\cat{C})$ provide a visual representation of
equations between these generalized elements. The word ``diagram'' is used here in the
ordinary sense of category theory. Thus, a \emph{diagram} in $\El(\cat{C})$ is a
functor $D : \cat{J} \to \El(\cat{C})$ from a small category $\cat{J}$, the
\emph{shape} or \emph{indexing category} of the diagram, to $\El(\cat{C})$. When the indexing category of a diagram is the discrete exterior calculus, we call that diagram a Decapode.
We choose $\cat{C}$ to be the category of smooth vector bundles over $M$, with objects of $\cat{C}$ differential forms on $M$. We write $\Omega^k(M)$ to represent the differential form of degree $k$ on $M$.

An immediate advantage of encoding equations as diagrams is that we can begin to interpret morphisms between such diagrams. A morphism of diagrams is given in Fig.~\ref{fig:diffusion_bcs}, which applies boundary conditions to the diffusion diagram given in Fig. \ref{fig:diff_theory_dec}.

\begin{figure}
    \centering
    \adjustbox{scale=0.9,center}{
\begin{tikzcd}
  {\Omega^0(\partial M)} & {\Omega^0(M)} && {\Omega^0_t(M)} & {\Omega^0_t(M)} & {\tilde\Omega^3_t(M)} \\
  & {\Omega_t^0(\partial M)} && {\Omega^1_t(M)} && {\tilde\Omega^2_t(M)}
  \arrow["{\partial_t}", from=1-4, to=1-5]
  \arrow["d"', from=1-4, to=2-4]
  \arrow["k\star"', from=2-4, to=2-6]
  \arrow["{\tilde d}"', from=2-6, to=1-6]
  \arrow["{\star^{-1}}"', from=1-6, to=1-5]
  \arrow["{\mathrm{res}_0}"', dashed, from=1-4, to=1-2]
  \arrow["{\mathrm{res}_{\partial M}}", dashed, from=1-4, to=2-2]
  \arrow["{\mathrm{res}_{{\partial M},0}}"', curve={height=18pt}, dashed, from=1-4, to=1-1]
  \arrow["{\mathrm{res}_{\partial M}}"', from=1-2, to=1-1]
  \arrow["{\mathrm{res}_0}", from=2-2, to=1-1]
\end{tikzcd}
    }
    \caption{Diagram morphism which both shows the underlying physical system and spatiotemporal boundary conditions. Solutions to this system are solutions to the related extension-lifting problem. ~\cite{patterson_diagrammatic_2023} (Example 4.10).}
    \label{fig:diffusion_bcs}
\end{figure}

Another feature of this diagrammatic syntax is that these diagrams are not spatially restricted. We can always project them onto different graphical layouts for application-specific visualization.

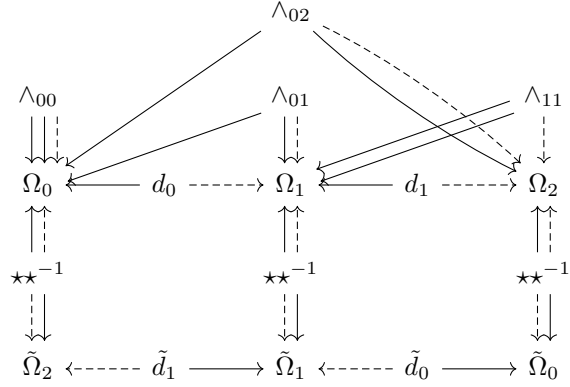
\begin{figure}[htbp]
\[\begin{tikzcd}
	&& {\wedge_{02}} \\
	{\wedge_{00}} && {\wedge_{01}} && {\wedge_{11}} \\
	{\Omega_0} & {d_0} & {\Omega_1} & {d_1} & {\Omega_2} \\
	{{\star}{\star^{-1}}} && {{\star} {\star^{-1}}} && {{\star} {\star^{-1}}} \\
	{\tilde{\Omega}_2} & {\tilde{d}_1} & {\tilde{\Omega}_1} & {\tilde{d}_0} & {\tilde{\Omega}_0}
	\arrow[shift left=1, from=4-1, to=3-1]
	\arrow[shift right=1, dashed, from=4-1, to=5-1]
	\arrow[shift right=1, dashed, from=4-1, to=3-1]
	\arrow[shift left=1, from=4-1, to=5-1]
	\arrow[shift right=1, dashed, from=4-5, to=3-5]
	\arrow[shift right=1, dashed, from=4-5, to=5-5]
	\arrow[shift left=1, from=4-5, to=5-5]
	\arrow[shift left=1, from=4-5, to=3-5]
	\arrow[shift right=1, dashed, from=4-3, to=5-3]
	\arrow[shift left=1, from=4-3, to=5-3]
	\arrow[shift right=1, dashed, from=4-3, to=3-3]
	\arrow[shift left=1, from=4-3, to=3-3]
	\arrow[shift right=1, from=2-5, to=3-3]
	\arrow[shift left=1, from=2-5, to=3-3]
	\arrow[dashed, from=2-5, to=3-5]
	\arrow[shift right=1, from=2-3, to=3-3]
	\arrow[shift left=1, from=2-3, to=3-1]
	\arrow[shift left=1, dashed, from=2-3, to=3-3]
	\arrow[shift left=3, dashed, from=2-1, to=3-1]
	\arrow[shift left=1, from=2-1, to=3-1]
	\arrow[shift right=1, from=2-1, to=3-1]
	\arrow[from=3-2, to=3-1]
	\arrow[dashed, from=3-2, to=3-3]
	\arrow[from=3-4, to=3-3]
	\arrow[dashed, from=3-4, to=3-5]
	\arrow[from=5-4, to=5-5]
	\arrow[dashed, from=5-4, to=5-3]
	\arrow[from=5-2, to=5-3]
	\arrow[dashed, from=5-2, to=5-1]
	\arrow[from=1-3, to=3-1]
	\arrow[curve={height=6pt}, from=1-3, to=3-5]
	\arrow[curve={height=-6pt}, dashed, from=1-3, to=3-5]
\end{tikzcd}\]
\caption{The hypergraph encoding of the de Rham complex as a presentation of a symmetric monoidal category. The objects  $\Omega_i$ are spaces of primal $i$-forms and $\tilde{\Omega}_i$ are the spaces of dual $i$-forms. All other objects denote operators in the discrete exterior calculus. Solid arrows point from operators to their inputs. Dashed arrows point from operators to their outputs. Wedge products $\wedge_{ij}$ are binary operators.}
\label{fig:deRahmSchema}
\end{figure}

The technique of $Hyp_\Sigma$ encoding ~\cite{10.1145/2933575.2935316} allows one to represent hypergraphs of physical computations as databases over the schema shown in Fig.~\ref{fig:deRahmSchema}. Databases over this schema can store the hypergraphs in the monoidal signature represented by the de Rham complex. In the Decapode schema, Fig. ~\ref{fig:decapode_schema_fig}, unary operators are encoded by tables with two foreign keys and binary operators are encoded with tables containing three foreign keys. In this schema, the operators are segregated into their own tables. So there is a separate table for each operator. For example, the $d_0$ table stores all the occurrences of the exterior derivative from $\Omega_0 \to \Omega_1$.

We can take a hypergraph and automatically generate a Julia program that implements the computation. Each vertex becomes a variable and each operator becomes a function call. The ordering of the function calls is consistent with the data flow constraints encoded in the hypergraph as described below. The program is constructed using the metaprogramming features of Julia to convert expressions representing equations into executable programs.
The Julia differential equations libraries \cite{rackauckas2017differentialequations} are then invoked with this function to evaluate the vector field. Boundary conditions are enforced via the insertion of extra operators that apply a masking operation to force the state variables on the boundary to have the value specified by the boundary condition data.

To schedule the execution of the hypergraph into a totally ordered sequence of function calls, we use a topological sorting algorithm to construct a total order on the set of variables and operations such that every variable is ordered after the operation that computes it and every operation is ordered after all the inputs to that operation.

The algorithm is given in pseudocode in Algorithm 1.
\begin{algorithm}
\begin{algorithmic}
\State{$\texttt{d} \gets$ the Decapode to be compiled}
\State{$\texttt{input\_numbers} \gets$ row numbers of the state variables in $Var$}
\State{$\texttt{visited} \gets $ array of $sizeof(Var)$ falses}
\State{$\texttt{visited[input\_numbers]} \gets$ true}
\State{$\texttt{consumed1} \gets $ array of $sizeof(Op1)$ falses}
\State{$\texttt{consumed2} \gets $ array of $sizeof(Op2)$ falses}
\State{$\texttt{consumed}\Sigma \gets $ array of $sizeof(\Sigma)$ falses}
\State{$\texttt{op\_order} \gets $ empty array}
\For{$i = 1 \to i =$ total size of the $Op_1$,$Op_2$, and $\Sigma$ tables}

\ForAll{rows $op \in Op1$}
\State{$\texttt{s} \gets d[op, :src]$ }
\If{$! \texttt{consumed1} [op] \textbf{ and } \texttt{visited} [s]$}
\State{$\texttt{operator} \gets d[op, :op1]$}
\State{$\texttt{t} \gets d[op, :tgt]$}
\If{$\texttt{operator} == \partial_t$}
  \State{continue}
\EndIf
\State{$\texttt{consumed1}[op] \gets$ true}
\State{$\texttt{visited}[op] \gets$ true}
\State{$\texttt{sname} \gets d[s, :name]$}
\State{$\texttt{tname} \gets d[t, :name]$}
\State{$\texttt{c} \gets UnaryCall(operator, sname, tname)$}
\State{$push!(\texttt{op\_order}, \texttt{c})$}
\EndIf
\EndFor

\ForAll{rows $op \in Op2$}
\State{$\texttt{arg1} \gets d[op, :proj1]$ }
\State{$\texttt{arg2} \gets d[op, :proj2]$ }
\If{$! \texttt{consumed1} [op] \textbf{ and } \texttt{visited} [arg1] \textbf{ and } \texttt{visited} [arg2]$}
\State{$\texttt{r} \gets d[op, :res]$}
\State{$\texttt{a1name} \gets d[arg1, :name]$}
\State{$\texttt{a2name} \gets d[arg2, :name]$}
\State{$\texttt{rname} \gets d[r, :name]$}
\State{$\texttt{operator} \gets d[op, :op2]$}
\State{$\texttt{consumed2} [op] \gets$ true}
\State{$\texttt{visited} [r] \gets$ true}
\State{$\texttt{sname} \gets d[s, :name]$}
\State{$\texttt{tname} \gets d[t, :name]$}
\State{$\texttt{c} \gets BinaryCall(operator, sname, tname)$}
\State{$push!(\texttt{op\_order}, \texttt{c})$}
\EndIf
\EndFor

\label{compile}

\algstore{compilealg}
\end{algorithmic}
\caption{The Decapode compilation algorithm}
\end{algorithm}

\begin{algorithm}
\begin{algorithmic}
\algrestore{compilealg}

\ForAll{rows $op \in \Sigma$}
\State{$args \gets $ rows in $d[:summand]$ incident to row $op$ of $d[:summation]$}
\If{$!consumed\Sigma[op] \textbf{ and } $all of $visited[args])$ is true}
\State{$r \gets d[op, :sum]$}
\State{$argnames \gets d[args, :name]$}
\State{$rname \gets d[r, :name]$}
\State{$operator \gets :+$}
\State{$consumed\Sigma[op] \gets$ true}
\State{$visited[r] \gets$ true}
\State{$sname \gets d[s, :name]$}
\State{$tname \gets d[t, :name]$}
\State{$c \gets VarargsCall(operator, argnames, rname)$}
\State{$push!(op\_order, c)$}
\EndIf
\EndFor

\EndFor

\State{$assigns \gets map(Expr, op\_order)$}
\State{$ret \gets :(return)$}
\State{$ret.args \gets d[d[:,:incl], :name]$}
\State{$f \gets $ new function with inputs $(du, u, p, t)$ wherein:}
\State{``````}
\State{Assign state variables to local variables}
\State{Unwrap $assigns$}
\State{Assign tangent variables to local variables}
\State{Return tangent variables}
\State{"""}

\hspace{-6.7mm}\Return{$f$}

\end{algorithmic}
\addtocounter{algorithm}{-1}
\caption{Decapode compilation algorithm, continued}
\end{algorithm}

\subsubsection{Requirements for successful simulator generation}

The algorithm above requires the following four conditions for successful generation of a simulator. We note that meeting these conditions does not guarantee that the generated simulator will be stable.

\begin{enumerate}
\item The directed bipartite graph encoding of the hypergraph must be acyclic. This guarantees that data dependencies are consistently ordered and that a topological sort is possible.

\item Every state variable must have exactly one corresponding time derivative to compute the update rule.

\item Every latent (non-state) variable must have at least one incoming edge to provide a value for that variable at each time step.

\item If a variable has more than one incoming edge, both edges must carry the same value. The software could be extended to relax this restriction by computing both edges and asserting the equality of the computed values. Alternatively, the software could be extended by taking the average of those computed values, so as to reduce error in the solution. This is a promising area for future work.
\end{enumerate}

Some physically valid diagrammatic equations cannot be solved with this computational approach. For example, if a variable is defined only implicitly as the variable $x$ for which $f(x) = y$, then the algorithm above will not successfully generate a solver. A simple example of this is given by an awkward presentation of the diffusion equation.  

\begin{figure}[hbtp]
    \centering
    \begin{subfigure}[b]{0.4\textwidth}
        \centering
        \begin{equation*}
        \begin{split}
            {\partial_t}C &= \Div \phi \\
            \phi &= k \Grad C
        \end{split}
        \end{equation*}
        \[\begin{tikzcd}
        	{C:\R_\Gamma} && {\phi:\R_\Gamma^2} \\
        	{\dot C:\R_\Gamma}
        	\arrow["k\nabla", from=1-1, to=1-3]
        	\arrow["{\partial_t}"', from=1-1, to=2-1]
        	\arrow["\nabla\cdot", from=1-3, to=2-1]
        \end{tikzcd}\]
        \caption{}
        \label{fig:comp_diffusion}
    \end{subfigure}
    \hspace{0.1\textwidth}
    \begin{subfigure}[b]{0.4\textwidth}
    \centering
    \begin{equation*}
    \begin{split}
        {\partial_t}C &= \Div \phi \\
        k^{-1}\phi &= \Grad C
    \end{split}
    \end{equation*}
        \[\begin{tikzcd}
        	{C:\R_\Gamma} & {\bullet:\R^2_\Gamma} \\
        	{\dot C:\R_\Gamma} & {\phi:\R_\Gamma^2}
        	\arrow["{\partial_t}"', from=1-1, to=2-1]
        	\arrow["\nabla\cdot", from=2-2, to=2-1]
        	\arrow["\nabla", from=1-1, to=1-2]
        	\arrow["k^{-1}"', from=2-2, to=1-2]
        \end{tikzcd}\]
    \caption{}
    \label{fig:uncomp_diffusion}
    \end{subfigure}
    \caption{These two equivalent descriptions of thermal diffusion highlight the structural requirements for a Decapode to be converted into a directed wiring diagram. In one presentation $\phi$ is computed from the $\nabla C$ value (\ref{fig:comp_diffusion}) and all computations can be scheduled. In the alternative presentation $\phi$ is only defined implicitly as $k^{-1}\phi = \nabla C$ (\ref{fig:uncomp_diffusion}). The algorithm for scheduling a hypergraph requires that there be an explicit definition of $\phi$. Note that the algebraic property of the equations is reduced to a combinatorial property of the graph.}
    \label{fig:comp_examples}
\end{figure}

Figure~\ref{fig:comp_examples} shows two equivalent equations, one which has the appropriate structure and one which does not. Unlike in Figure~\ref{fig:comp_diffusion}, the system of equations in Figure~\ref{fig:uncomp_diffusion} does not satisfy the structural requirements because $\phi$ does not appear independently on one side of an equation. The presentation of diffusion shown in Figure~\ref{fig:uncomp_diffusion} can be  rewritten into that of Figure~\ref{fig:comp_diffusion} using the fact that $k$ is an invertible linear map (as long as the diffusion coefficient is nonzero everywhere). The theory of diagrammatic equations allows mathematicians to reason about these relationships between presentations of physical equations, which can then be used in software to relate numerical simulations. In this case, two equivalent mathematical models are shown that lead to different computational methods. It should be noted that this is by no means the only method of solving Decapodes, but it is one which provides accurate numerical results for certain problems.

While the theory of diagrammatic equations handles boundary conditions using an elegant formalism based on morphisms in the category of diagrams~\cite{patterson_diagrammatic_2023}, the Decapodes software system uses a more traditional approach. In order to enforce boundary conditions, the relevant value is restricted to the boundary condition after all other computations are performed. For example, the boundary condition $${\partial_t} T(x,y) = 0 \qquad\text{where}\qquad (x,y)\in \partial X$$
would be enforced by setting the value of $\dot{T}$ to zero on the boundary points after $\dot{T}$ has been computed with a vector masking operation. This ensures that the boundary condition is satisfied, at the expense of spatial continuity in the value or satisfaction of previous equalities. While the numerical properties of this technique have not yet been investigated, this has proved effective for computing quantities in the benchmarks discussed in section \ref{sec:experiments}.

\subsection{DEC Formulation of Problems}
\paragraph{Partial Implementation of Compressible Navier-Stokes}
Using the exterior calculus instead of vector calculus to formulate systems of equations in multiphysics
has two major benefits: it
generalizes the equations from Euclidean space to an arbitrary Riemannian
manifold, and it enables a direct translation to the discrete exterior calculus
and thus to the Decapodes software for numerical simulations. Although fluid mechanics has been studied using the exterior calculus~\cite[\S
8.2]{abraham1988} and its discrete counterpart~\cite{wilson2011,Mohamed_2016},
particularly in the case of incompressible fluids, the general form of the
Navier-Stokes equations for compressible, viscous fluids seems not to be treated
from this viewpoint in the literature.
We propose here an exterior calculus formulation of the compressible Navier-Stokes momentum equation in eq. \ref{eq:ns_ec_momentum}. However, we note here that we do not include the term due to viscous dissipation in the energy balance eq. \ref{eq:ns_ec_energy_balance}, which the updated momentum equation requires to accurately capture compressible flow. A full treatment of the compressible Navier-Stokes equations in the exterior calculus that uses an updated energy balance equation is grounds for further work.
\begin{figure}
    \centering
    \begin{subfigure}[b]{0.4\textwidth}
        \begin{align}
  \frac{\Dif \rho}{\Dif t} &= -\rho \Div \vec{u} \\
  \rho \frac{\Dif \vec{u}}{\Dif t} &=
    - \Grad p + \mu \nabla^2 \vec{u} + \frac{1}{3} \mu \Grad(\Div \vec{u}) \\
  \rho \frac{\Dif e}{\Dif t} &= -p \Div \vec{u}
\end{align}
        \caption{Navier-Stokes as expressed in vector calculus}
        \label{fig:ns_vc}
    \end{subfigure}
    \hspace{0.1\textwidth}
    \begin{subfigure}[b]{0.4\textwidth}
        \begin{align}
  \Dif_t \tilde e &= (\rho \tilde e - \tilde p) \delta u\\
  \rho \Dif_t u &= -\dif p - \mu \Delta u + \frac{1}{3} \mu \dif \delta u \label{eq:ns_ec_momentum} \\
  \rho \Dif_t \tilde e &= (\rho \tilde e - \tilde p) \delta u\label{eq:ns_ec_energy_balance}
\end{align}
        \caption{Navier-Stokes as expressed in exterior calculus}
        \label{fig:ns_dec}
    \end{subfigure}
    \caption{While the operators used differ between the vector calculus and exterior calculus expressions of the Navier-Stokes equations, the structure of the equations remain similar. There are some modeling choices that must be made going from vector to exterior calculus, but the process is mostly mechanical.}
    \label{fig:ns_expr}
\end{figure}
Recall that the Navier-Stokes equations in convective form are of the form given in Figure \ref{fig:ns_vc}, where
$\frac{\Dif}{\Dif t} := \frac{\partial}{\partial t} + \vec{u} \cdot \nabla$ is
the material derivative, $\rho$ is the mass density, $\vec u$ is the velocity
vector field, $p$ is the pressure, and $e$ is the specific internal energy.

We propose the generalized equations in Figure \ref{fig:ns_dec},
where $u := \vec{u}^\flat$ is the 1-form proxy for the vector field $\vec{u}$
and $\tilde \rho$ and $\tilde e$ are the twisted $n$-forms obtained by Hodge
duality from the 0-forms $\rho$ and $e$. The material derivative
$\Dif_t := \partial_t + \nabla_{\vec u}$ on a Riemannian manifold is defined
using the covariant derivative $\nabla_{\vec u}$~\cite{boyland2001}, which in
the special cases of interest has expressions involving the Lie derivative
$\mathcal{L}_{\vec u}$. When the manifold is Euclidean space $\R^3$
with its standard metric and coordinate system, the above equations reduce to
the classical Navier-Stokes equations.

\section{Experimental results}
\label{sec:experiments}

A variety of numerical experiments are performed for comparing various aspects of the performance of Decapodes to a state-of-the-art solver, SU2. The following sections detail the characteristics of SU2, along with the selected modeling problems. 

\subsection{Benchmarks}
We use standard challenge problems to benchmark the Decapodes tool against, using conjugate heat transfer as our primary target. We compare the results of Decapodes against the state of the art tool SU2~\cite{EconomonSU2}, described below. These benchmark scripts are made available online\footnote{\url{https://github.com/AlgebraicJulia/DECAPODES-Benchmarks/tree/main/simulations}}. We note that Keyes et al.~\cite{keyes_multiphysics_2013} offer several classes of coupling that may constitute a "multiphysics" system, such as those that couple along an interface. For the purposes of demonstrating this software framework, studying the current style of term-coupling will suffice.

It should be noted that the current method of compiling a Decapode into a program via topologic sort results in an explicit time-stepping method. This means that in order to evaluate these physics on a large mesh, we need to use very small timesteps to ensure stability (much smaller than those used by the library we are comparing against). Since the use of implicit and more advanced methods would be a significant but straightforward engineering effort, the high runtime of our tool with respect to SU2 is expected and manageable; accuracy of simulations is of primary interest.

\subsubsection{SU2 Simulation Framework}
In order to test the performance of Decapodes, we compare against an existing state-of-the-art simulation framework, namely SU2. SU2~\cite{EconomonSU2} is an open source library for the simulation of multiphysics problems comprised of partial differential equations on arbitrary unstructured meshes. Of primary interest for this effort is the SU2\_CFD (computational fluid dynamics) executable, designed for the solution of complex scalar transport PDEs such as the Euler, Navier-Stokes, Poisson, and Heat equation among others. SU2 is executable in both serial and parallel modes via MPI and mesh partitioning developed around parMETIS~\cite{EconomonSU2}. The SU2\_CFD solver likewise is able to operate as a Finite Volume or Finite Element solver. Additionally, various time integration schemes, both implicit and explicit are available for use. Simulations are defined via input specification files along with meshes of the computational domain. For this study, the \Deca tool utilizing a DEC backend is compared to the SU2\_CFD FVM solver, testing various coupled transport phenomena. 

\subsubsection{Conjugate Heat Transfer}
SU2 features support for multiphysics simulations comprising of compositions of their fluid mechanics, structural mechanics, and heat transfer solvers. Tutorial examples include generation and solution of fluid-structure-interaction simulations via the coupling of the Navier-Stokes equations with the structural mechanics solver, and conjugate heat transfer (CHT) via the coupled Navier-Stokes and Energy Equations over the fluid domain, to the Heat Equation over solid regions. 

The purpose of this benchmark is to compare the compositionality of Decapodes to generating multiphysics simulations in SU2. For SU2, a special interface for multiphysics \texttt{SOLVER= MULTIPHYSICS} is used to allow for the coupling of different zones of computation, each with their own corresponding problem definitions and subconfiguration and denoted via the \texttt{CONFIG\_LIST} command, through their multizone simulation framework. For the conjugate heat transfer problem, coupling is done through designated interfaces at the boundaries of each subdomain via the \texttt{MARKER\_CHT\_INTERFACE}. 
Note that the SU2 has conjugate heat transfer implemented as a hard-coded, predefined physics with a corresponding solver implementation provided by the SU2 developers. In the Decapodes.jl implementation, the user is the one who defines the CHT model, and the software generates the solver implementation. This represents a radically different approach to multiphysics simulations, emphasizing a flexibility in creating novel models.

\begin{figure}[htbp]
    \centering
    \includegraphics[width=0.5\textwidth]{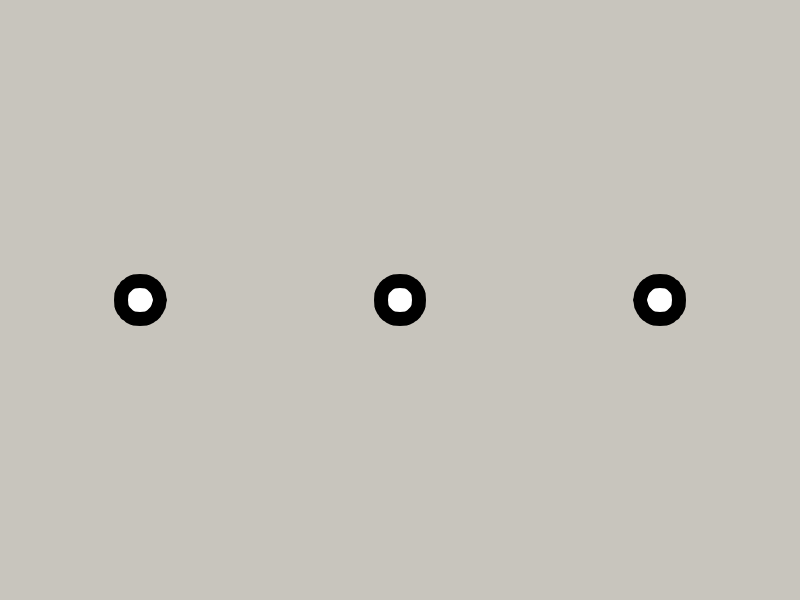}
    \caption{CHT Schematic. Black regions represent the solid domain cylinders. Gray regions represents the fluid domain.}
    \label{fig:cht_schematic}
\end{figure}

\paragraph{Domain Region}
This test problem involves a fully coupled multiphysics example of conjugate heat transfer, featuring non-isothermal flow around a set of heated solid cylinders. This test uses the incompressible Navier-Stokes formulation of eqn. \ref{eq:incomp_ns}, resolving the energy equation for temperature in the fluid field, while the heat equation is resolved in the solid cylinders as defined by eqn. \ref{eq:heat_eqn}.
\begin{equation}\label{eq:heat_eqn}
\partial_t U - \nabla \cdot \bar{F}^v(U, \nabla U) - S = 0
\end{equation}
\begin{equation}
    U = \{\rho_s c_{p,s} T \}, \: \bar{F}^v(U, \nabla U) = \kappa_s \nabla T,
\end{equation}
with $S$ a source.
A schematic of the computational geometry is provided as Figure \ref{fig:cht_schematic}. Fluid enters the domain from the left boundary at a velocity and temperature of $U_{inlet} = 3.40297 [\mathtt{\frac{m}{s}}]$ and $T_{inlet} = 288.15 [\mathtt{K}]$. Each solid cylinder has a Dirichlet boundary condition of $T_{core} = 350 [\mathtt{K}]$. 

\paragraph{Physics}

\begin{figure}[htbp]
    \centering
    \includegraphics[width=0.5\columnwidth]{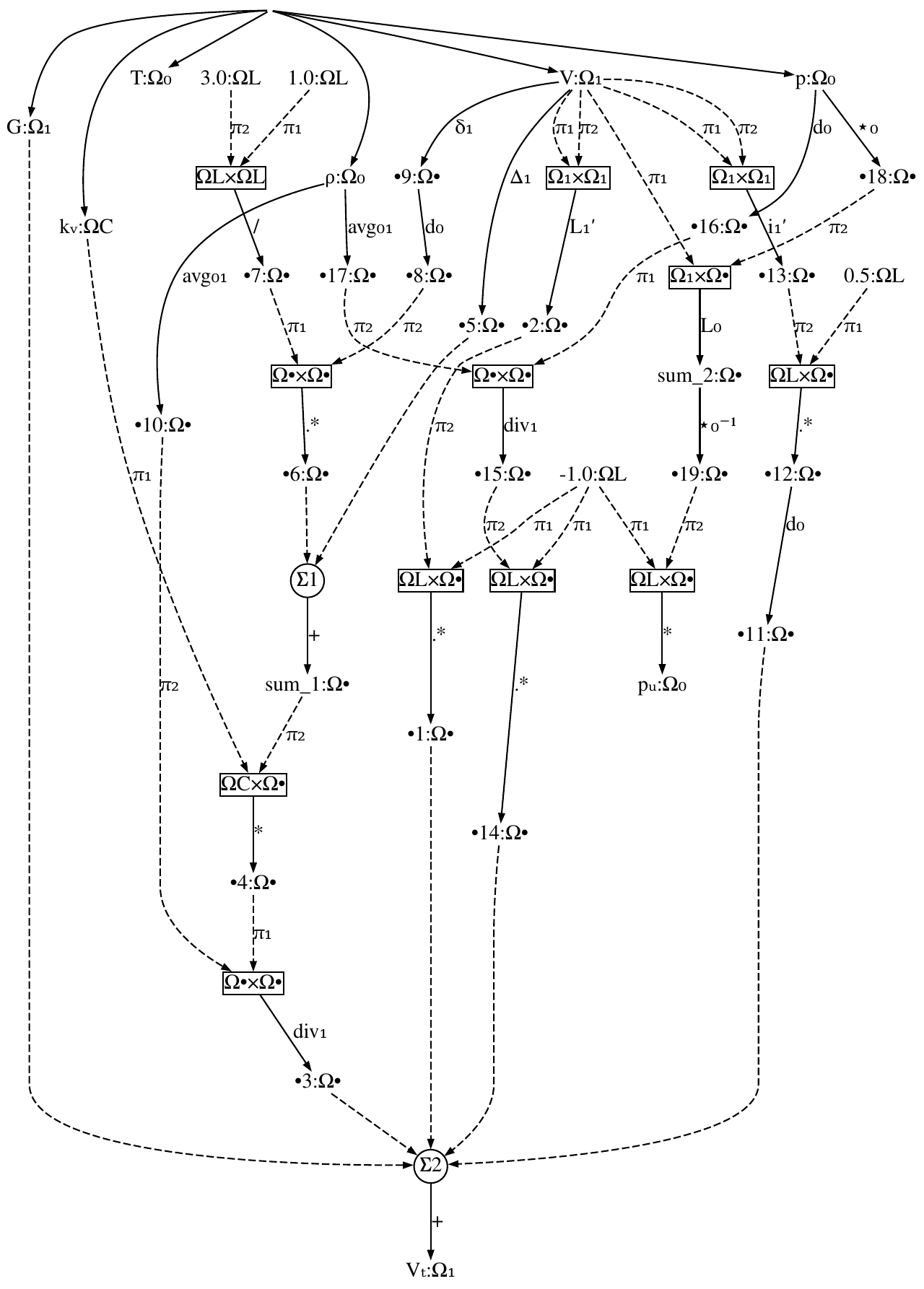}
    \caption{The CHT Navier-Stokes equations Decapode. Note that the State variables are denoted as those which are pointed to by some invisible node. Our layout algorithm naturally places such variables near the top of the figure. One may observe that these formal Decapodes diagrams are suggestive of the visualizations of dependency graphs and dataflow programs.}
    \label{fig:cht_navierstokes}
\end{figure}

\begin{figure}[htbp]
    \centering
    \begin{subfigure}[t]{0.2\textwidth}
    \includegraphics[width=\textwidth]{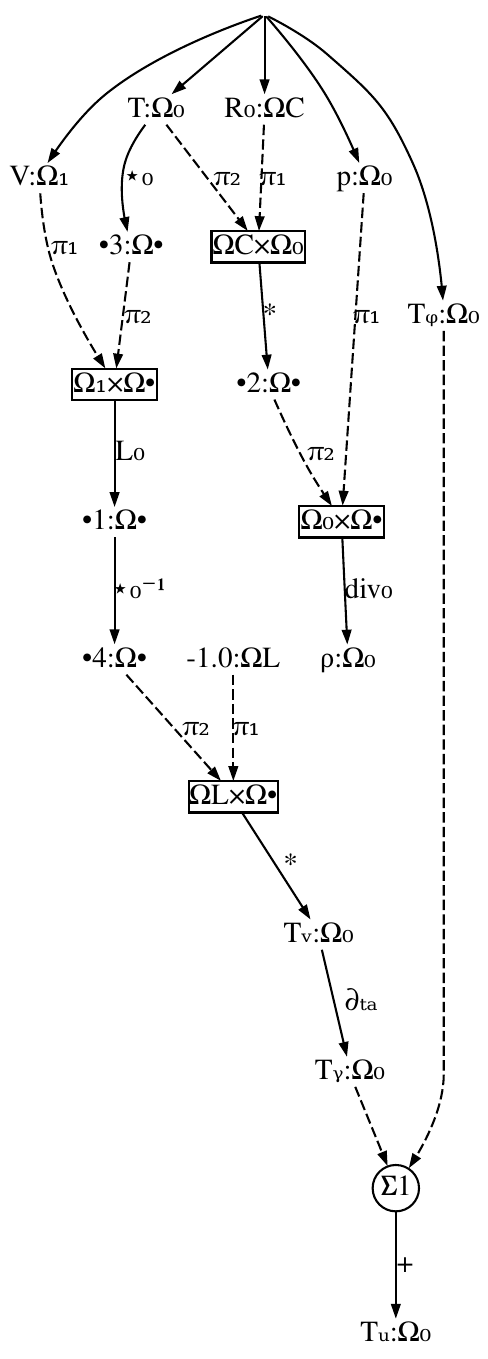}
    \caption{}
    \label{fig:CHT_Energy_Decapode}
    \end{subfigure}
    \hspace{1cm}
    \begin{subfigure}[t]{0.1\textwidth}
    \includegraphics[width=\textwidth]{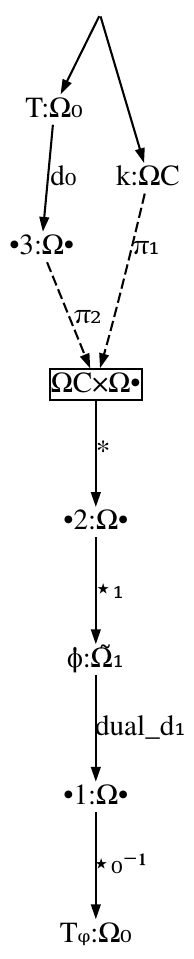}
    \caption{}
    \end{subfigure}
    \hspace{1cm}
    \begin{subfigure}[t]{0.3\textwidth}
    \includegraphics[width=\textwidth]{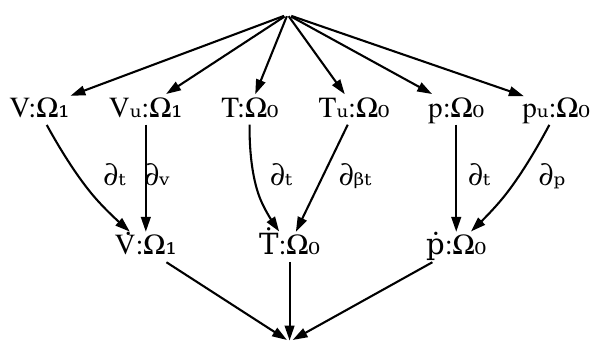}
    \caption{}
    \end{subfigure}
\caption{The CHT energy balance equations Decapode, diffusion equation Decapode, and boundary equations Decapode, respectively. We denote tangent variables here by drawing arrows from such nodes to a point at the bottom of the figure.}
    \label{fig:cht_energy_diffusion_boundaries}
\end{figure}

\begin{figure}[htbp]
    \centering
    \includegraphics[width=0.5\textwidth]{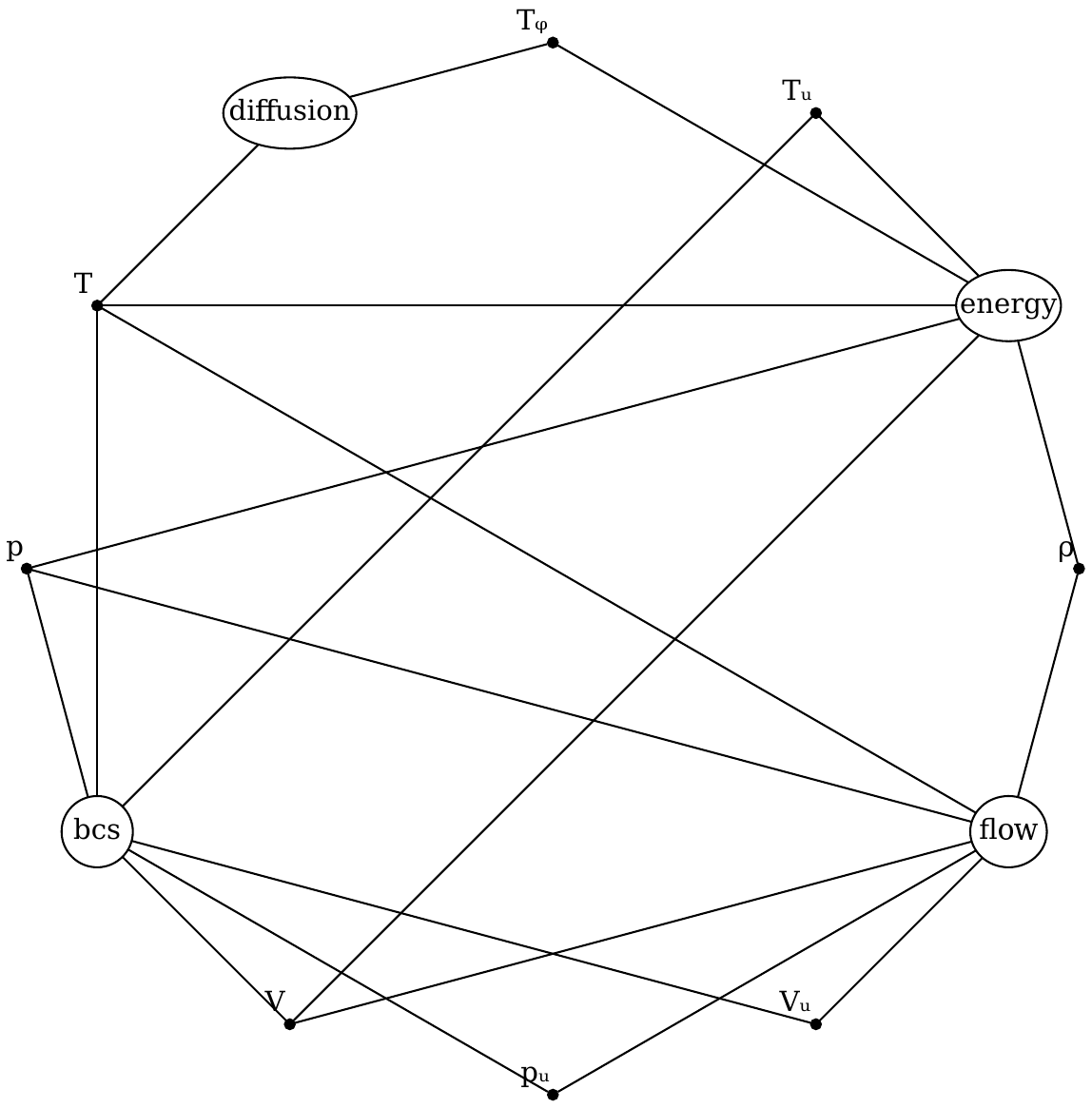}
    \caption{Composition pattern for the CHT benchmark problem. This provides a composition pattern for the Decapodes describing thermal diffusion, the energy balance equations, the fluid flow, and the problem boundary conditions. This visualization is algorithmically generated from the description of the composition pattern in code, and is layed out using Graphviz's ``circo" algorithm.}
    \label{fig:cht_pipe_uwd}
\end{figure}

\begin{figure}[htbp]
    \centering
    \begin{tiny}
    \includegraphics[width=0.6\columnwidth]{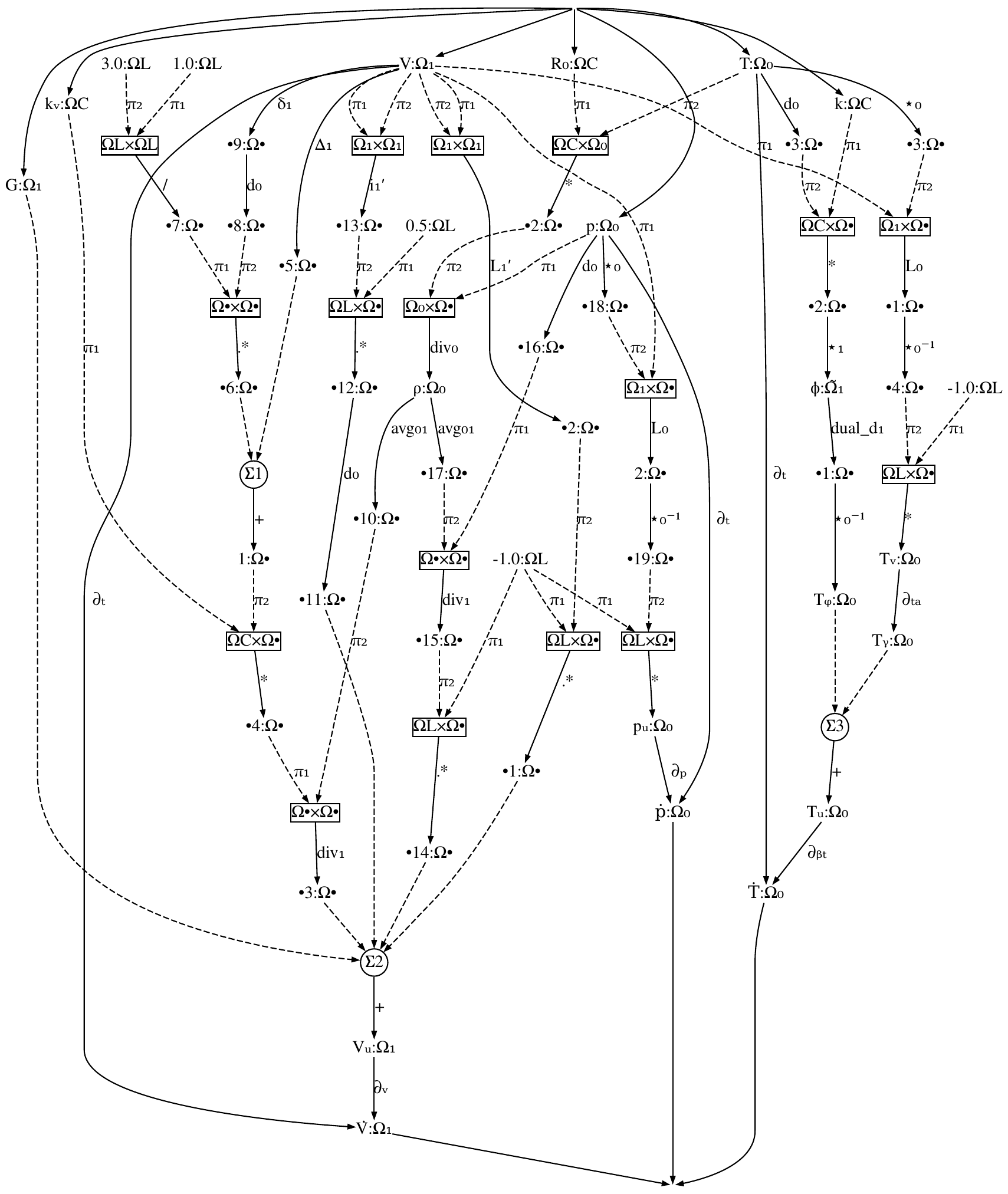}
    \caption{The final composed CHT Decapode.}
    \label{fig:cht_heatxfer_neato}
    \end{tiny}
\end{figure}

For the Decapodes implementation of this problem, we follow the DEC derivation of the 
incompressible Navier-Stokes equations
derived by Patterson et. al.~\cite{patterson_diagrammatic_2023}, except instead of using internal energy as a state variable, we use pressure, temperature, and velocity as our state variables, as shown in Fig.~\ref{fig:cht_navierstokes}. The Decapodes which encode thermal energy balance, diffusion, and boundary condition physics are given as in Fig. \ref{fig:cht_energy_diffusion_boundaries}, and are composed at the operadic level with the Navier-Stokes equations as in Fig. \ref{fig:cht_pipe_uwd}. This results in the Decapode shown in Fig. \ref{fig:cht_heatxfer_neato}, which is then transformed into its computation graph. This resulting program is then simulated according to the initial conditions and domain described above. The boundary conditions enforced by Fig. \ref{fig:cht_energy_diffusion_boundaries} set the change in pressure ``inside'' the cylinders to 0, the change in velocity along the walls of the cylinders and the walls of the domain to be 0, and the change in temperature along the walls of the cylinders and the walls of the domain to be 0.

\paragraph{Computational Results}

\begin{figure}[htbp]
     \centering
     \begin{subfigure}[htbp]{0.49\textwidth}
         \centering
         \includegraphics[width=0.99\textwidth]{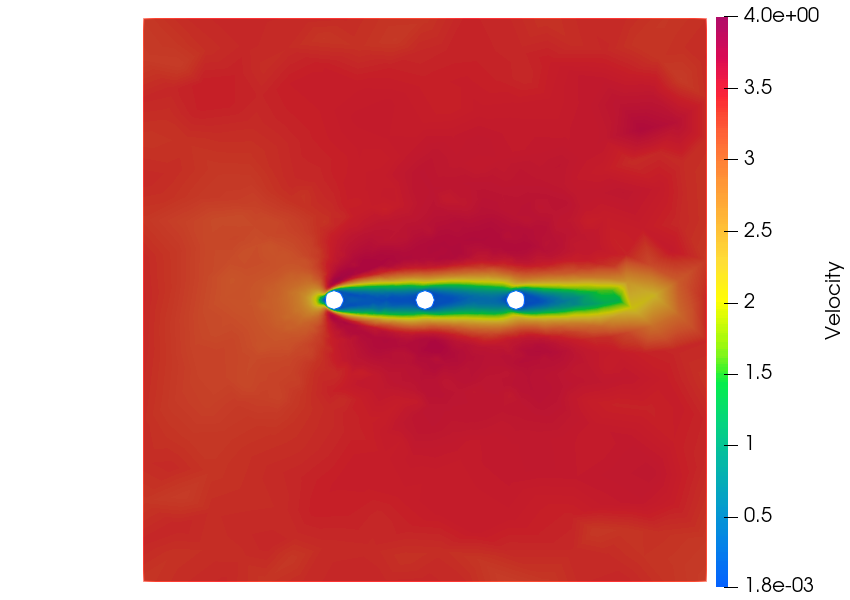}
         \caption{Velocity Field Computed with \Deca}
         \label{fig:cht_vel_dec}
     \end{subfigure}
     \hfill
     \begin{subfigure}[htbp]{0.49\textwidth}
         \centering
         \includegraphics[width=0.99\textwidth]{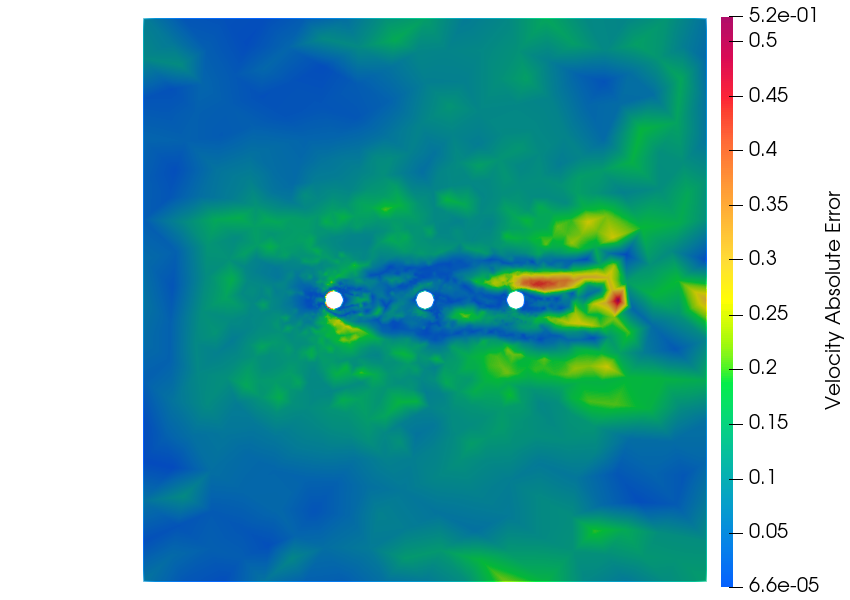}
         \caption{Absolute Error between Velocity Calculation of \Deca and SU2}
         \label{fig:cht_vel_error}
     \end{subfigure}
        \caption{CHT Benchmark Velocity Comparison}
        \label{fig:cht_vel_results}
\end{figure}

\begin{figure}[htbp]
     \centering
     \begin{subfigure}[htbp]{0.49\textwidth}
         \centering
         \includegraphics[width=0.99\textwidth]{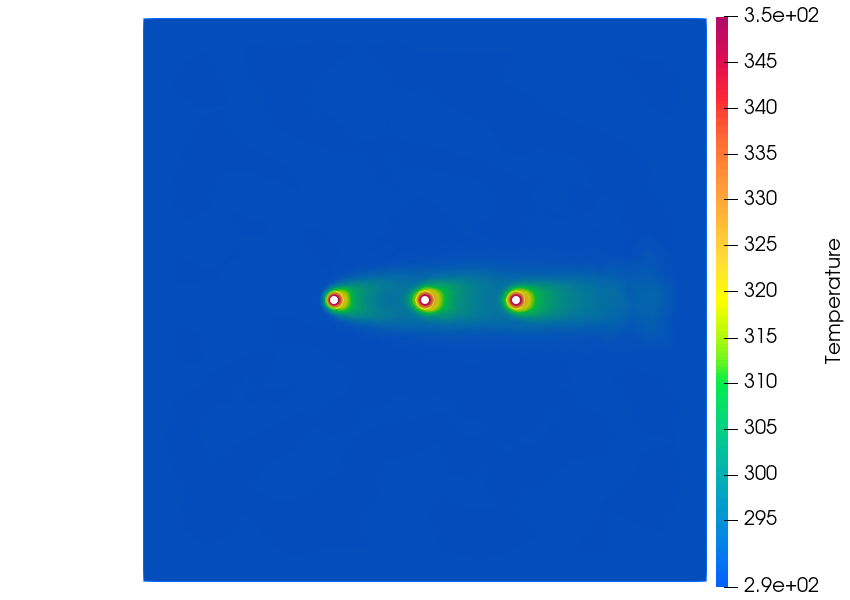}
         \caption{Temperature Field Computed with \Deca}
         \label{fig:cht_temp_dec}
     \end{subfigure}
     \hfill
     \begin{subfigure}[htbp]{0.49\textwidth}
         \centering
         \includegraphics[width=0.99\textwidth]{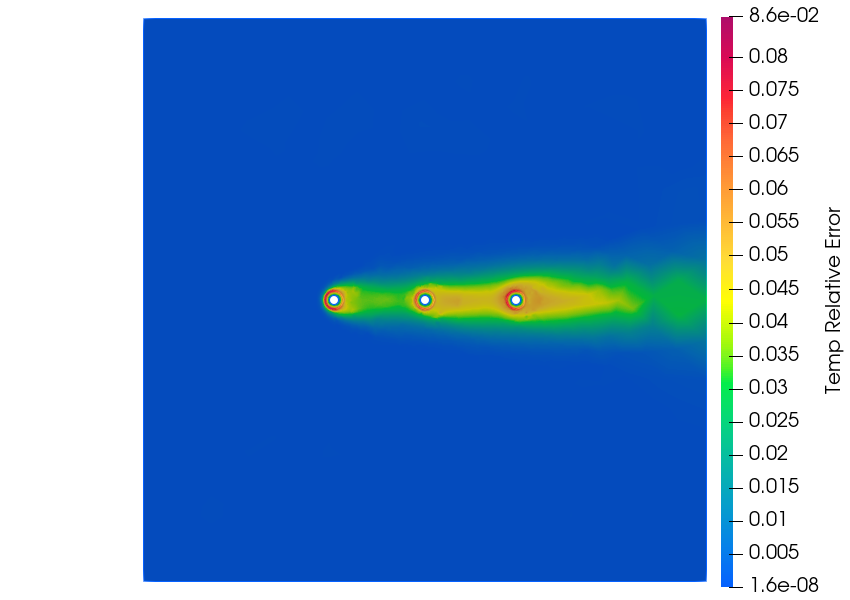}
         \caption{Relative Error between Temperature Calculation of \Deca and SU2}
         \label{fig:cht_temp_error}
     \end{subfigure}
        \caption{CHT Benchmark Temperature Comparison}
        \label{fig:cht_temp_results}
\end{figure}

\begin{figure}[htbp]
     \centering
     \begin{subfigure}[htbp]{0.4\textwidth}
         \centering
         \includegraphics[width=0.99\textwidth]{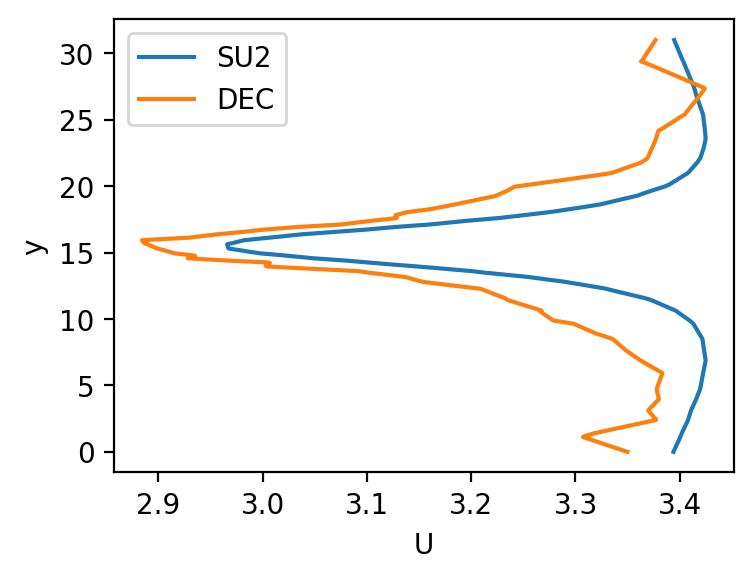}
         \caption{Velocity Profile along midline at x=-2.0 of domain.}
         \label{fig:cht_vel_profile_-2}
     \end{subfigure}
     \hfill
     \begin{subfigure}[htbp]{0.4\textwidth}
         \centering
         \includegraphics[width=0.99\textwidth]{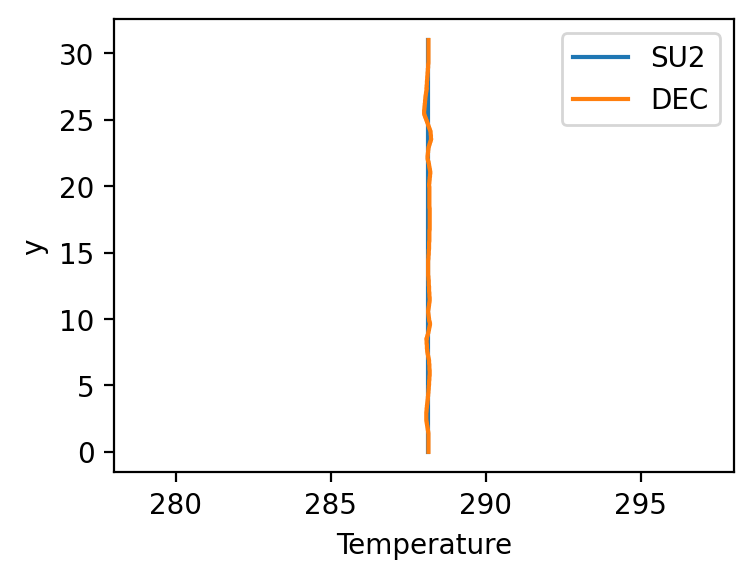}
         \caption{Temperature Profile along midline at x=-2.0 of domain.}
         \label{fig:cht_temp_profile_-2}
     \end{subfigure}
     
     \begin{subfigure}[htbp]{0.4\textwidth}
         \centering
         \includegraphics[width=0.99\textwidth]{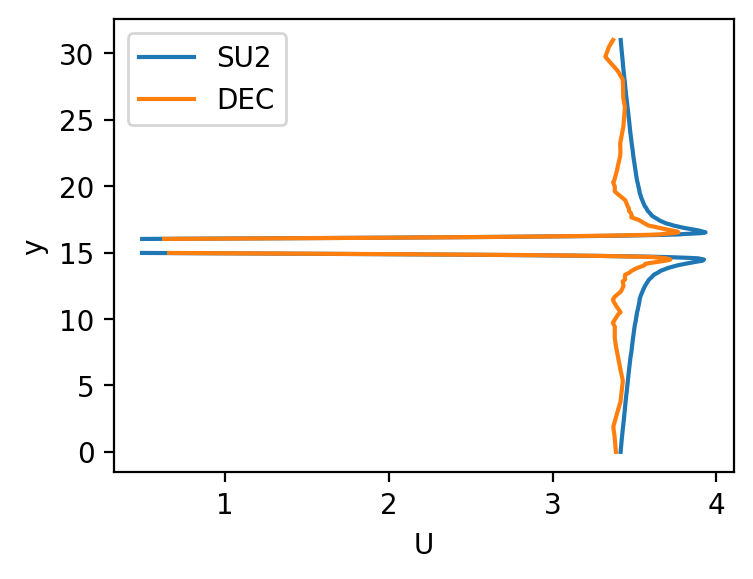}
         \caption{Velocity Profile along midline at x=0.5 of domain.}
         \label{fig:cht_vel_profile_0-5}
     \end{subfigure}
     \hfill
     \begin{subfigure}[htbp]{0.4\textwidth}
         \centering
         \includegraphics[width=0.99\textwidth]{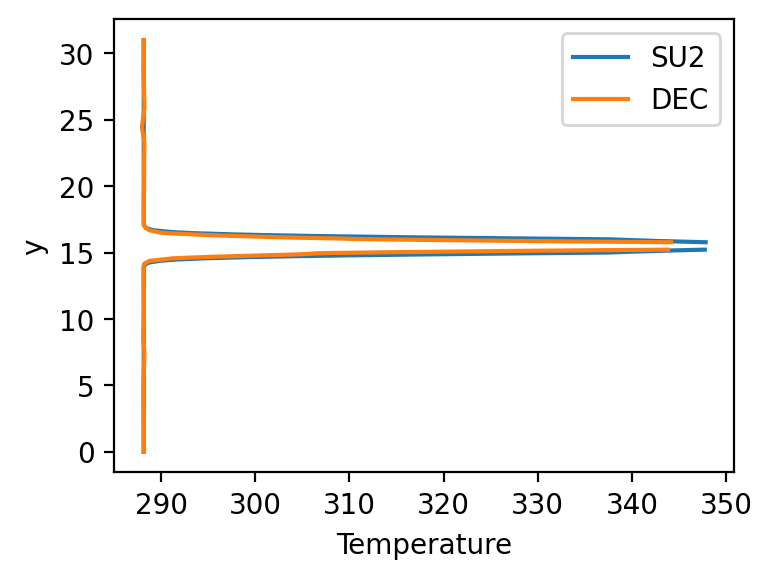}
         \caption{Temperature Profile along midline at x=0.5 of domain.}
         \label{fig:cht_temp_profile_0-5}
     \end{subfigure}
     
     \begin{subfigure}[htbp]{0.4\textwidth}
         \centering
         \includegraphics[width=0.99\textwidth]{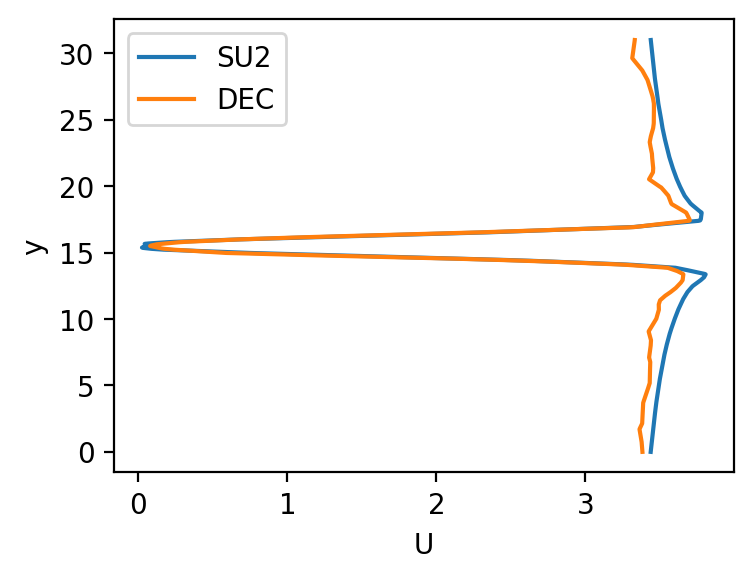}
         \caption{Velocity Profile along midline at x=3.0 of domain.}
         \label{fig:cht_vel_profile_3-0}
     \end{subfigure}
     \hfill
     \begin{subfigure}[htbp]{0.4\textwidth}
         \centering
         \includegraphics[width=0.99\textwidth]{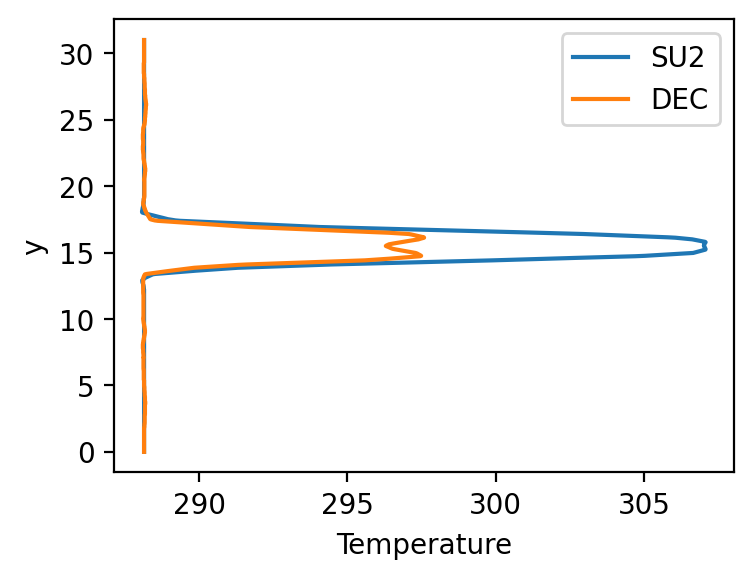}
         \caption{Temperature Profile along midline at x=3.0 of domain.}
         \label{fig:cht_temp_profile_3-0}
     \end{subfigure}
        \caption{Conjugate Heat Transfer Benchmark Comparison}
        \label{fig:cht_results_profile}
\end{figure}

For this benchmark, the simulation results from \Deca are compared directly to the SU2 results. The simulations are marched for a total time of $t = 5 \mathrm{[s]}$. Time marching is performed in SU2 via an implicit solver method to allow for larger timesteps. At the time of this writing, only explicit solution methods are available for the Decapodes formulation, with $\mathtt{Tsit5}()$ being used from the \texttt{DifferentialEquations.jl} library. The qualitative behavior of the velocity and temperature fields are examined in Figures \ref{fig:cht_vel_results} and \ref{fig:cht_temp_results} respectively. Decapodes is generally able to reconstruct the physical behavior in agreement with the results calculated via SU2. Both frameworks feature the creation of a wake in the flow behind each cylinder, and likewise cooling of the heated cylinders and subsequent heating of the fluid in response. Likewise as visible in Figures \ref{fig:cht_vel_error}, the absolute error between the velocity fields of Decapodes and SU2, and \ref{fig:cht_temp_error}, the relative error between the temperature fields calculated by Decapodes and SU2, remain generally low. Complete absolute and relative error for velocity, temperature, and density are given Table~\ref{tab:cht_dens_benchmarks}. The discrepancies between the velocity fields of SU2 and Decapodes likely comes from different initial conditions. SU2 performs a smoothing procedure to begin the simulation without
sharp initial conditions, whereas Decapodes has a pressure wave from the beginning which continues advecting through the fluid body during the simulation, causing variation in the velocity field.

More granular comparisons between the Decapodes and SU2 results are carried out via vertical profiles at locations $x = -2.0, 0.5, 3.0$ and presented in Figure \ref{fig:cht_results_profile} for the velocity and temperature profiles. Again, good agreement is visible between the 
Decapodes and SU2 results with the relative difference in the velocity profiles on the order of $5\%$ with the temperature values remaining within $\pm 3.5\%$ in the main frontal area of the flow, with a higher discrepancy in the wake region. This higher discrepancy is also likely due to varying initial conditions, where SU2 begins with a more developed temperature field than the 
Decapodes simulation. The lack of a smoothing procedure, and the choice of error tolerance selected for the $\mathtt{Tsit5}()$ algorithm are the likely causes for the sharpness of the selected temperature profiles. The under-prediction of the Decapodes results for temperature exactly between adjacent cylinders - such as demonstrated in Fig.~\ref{fig:cht_temp_profile_3-0} - may be due to similar reasons. This phenomenon may lessen as resolution increases, but we will not investigate the numeric stability as resolution increases for such DEC methods, and leave it instead as future work.

\paragraph{Runtime Results}
Computation time per second of simulation time for this benchmark is summarized in table \ref{tab:DEC_benchmarks}.
We believe the difference in total simulation time is because 
\Deca currently only uses explicit time-stepping methods which require significantly smaller time-steps to remain stable on fine meshes. This also explains the longer individual time-steps of SU2, as it takes longer to compute the stable solution for each time-step, while Decapodes simply performs a single application of the operators.

\subsubsection{Buoyancy-Driven Fluid Flow}
This test problem is of a cavity buoyancy-driven flow, which features natural convection driven by temperature induced buoyancy effects, while under the effect of gravity, for standard air modeled as an ideal gas.

The SU2 simulation uses an incompressible Navier-Stokes formulation with variable density as defined by eqn. \ref{eq:incomp_ns}.
\begin{equation}\label{eq:incomp_ns}
    \frac{\partial V}{\partial t} + \nabla \cdot \bar{F}^c(V) - \nabla \cdot \bar{F}^v(V, \nabla V) - S = 0
\end{equation}
\begin{equation}
    U = \begin{bmatrix}
    \rho \\
    \rho \bar{v} \\
    \rho c_p T
    \end{bmatrix}, \:
    V = \begin{bmatrix}
    p \\
    \bar{v} \\
    T
    \end{bmatrix}, \:
    \bar{F}^c(V) = \begin{Bmatrix}
    \rho \bar{v} \\
    \rho \bar{v} \otimes \bar{v} + \overline{\overline{I}}p \\
    \rho E \bar{v} + p \bar{v}
    \end{Bmatrix}, \:
    \bar{F}^v(V, \nabla V) = \begin{Bmatrix}
    \cdot \\
    \overline{\overline{\tau}} \\
    \overline{\overline{\tau}} \bar{v} + \kappa \nabla T
    \end{Bmatrix}
\end{equation}
\begin{equation}
    \overline{\overline{\tau}} = \mu (\nabla \bar{v} + \nabla \bar{v}^T) - \mu \frac{2}{3} \overline{\overline{I}} (\nabla \cdot \bar{v}),
\end{equation}
with $S$ a source.
In this formulation, the conserved variables of $U$ correspond to mass continuity, momentum balance, and the energy equation. However, the working variables for the incompressible solver are actually as defined in $V$ being the pressure, velocity, and temperature fields respectively. Density variations are determined as a function of temperature $\rho = \rho(T)$ via the ideal gas equation $\rho = P_0 / RT$ evaluated at reference pressure $P_0$. Likewise $\overline{\overline{\tau}}$ corresponds to the viscous stress tensor. In this simulation the primary parameter affecting the solution is the Rayleigh number defined as \cite{Sockol2003MultigridSO}:
\begin{equation}
    \mathrm{Ra} = \frac{\Delta T \rho^2 g c_p L^3}{\mu \kappa}
\end{equation}
Here $\Delta T$ is the mean temperature difference between the hot and cold walls:
\begin{equation*}
    \Delta T = 2 (T_h - T_c)/(T_h + T_c)
\end{equation*}
$\kappa$ is an additional scaled thermal conductivity parameter defined by the Prandtl number $\mathrm{Pr}$, and the material thermal conductivity $k$.
\begin{equation*}
    \kappa = \frac{\mu c_p}{\mathrm{Pr}}, \: \mathrm{Pr} = \frac{\mu c_p}{k}
\end{equation*}
The Rayleigh number is varied by changing the thermal conductivity of the fluid and is tested at the value $\mathrm{Ra} = 1\mathrm{e}5$.

\begin{figure}[htbp]
    \centering
    \includegraphics[width=0.5\textwidth]{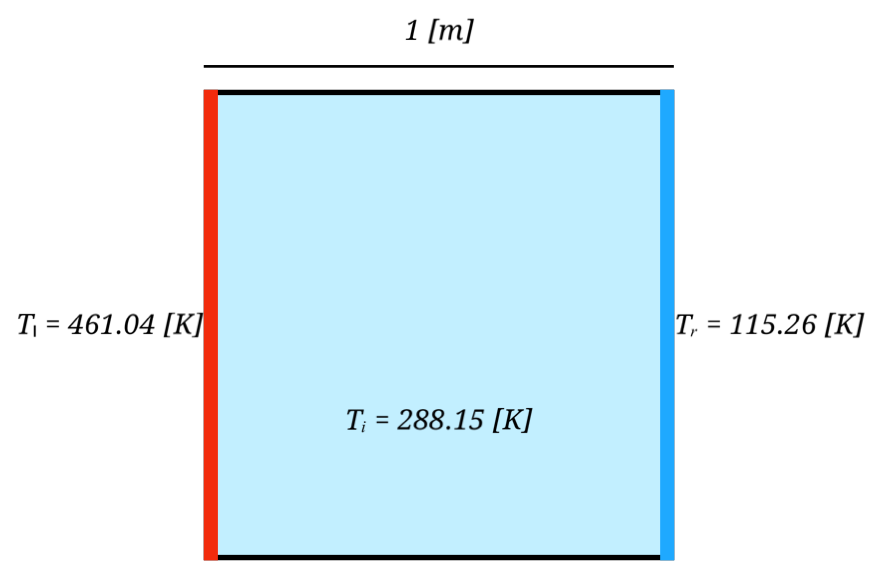}
    \caption{Schematic of the buoyancy-driven flow domain.}
    \label{fig:cavity_schematic}
\end{figure}

\paragraph{Domain Region}
The domain of computation is a $1 [\mathtt{m}]$ square cavity. Natural convection is driven by applying Dirichlet boundary conditions to the temperature field on the left and right walls of the cavity. The left boundary is set to $T_l = T_h = 461.04 \mathrm{[K]}$, while the right boundary is set to $T_r = T_c = 115.26 \mathrm{[K]}$. The top and bottom walls of the cavity are set to be adiabatic using Neumann boundary conditions. Lastly, the initial temperature in the cavity is set to $T_i = 288.15 \mathrm{[K]}$. A schematic of this domain is given in Fig. ~\ref{fig:cavity_schematic}. Buoyancy-driven fluid flow in a square cavity domain is a well known benchmark, with examples available in online documentation.\footnote{\url{https://reference.wolfram.com/language/PDEModels/tutorial/Multiphysics/ModelCollection/BuoyancyDrivenFlow.html}} 

\paragraph{Physics}

\begin{figure}[htbp]
    \centering
    \begin{tiny}
    \includegraphics[width=0.5\columnwidth]{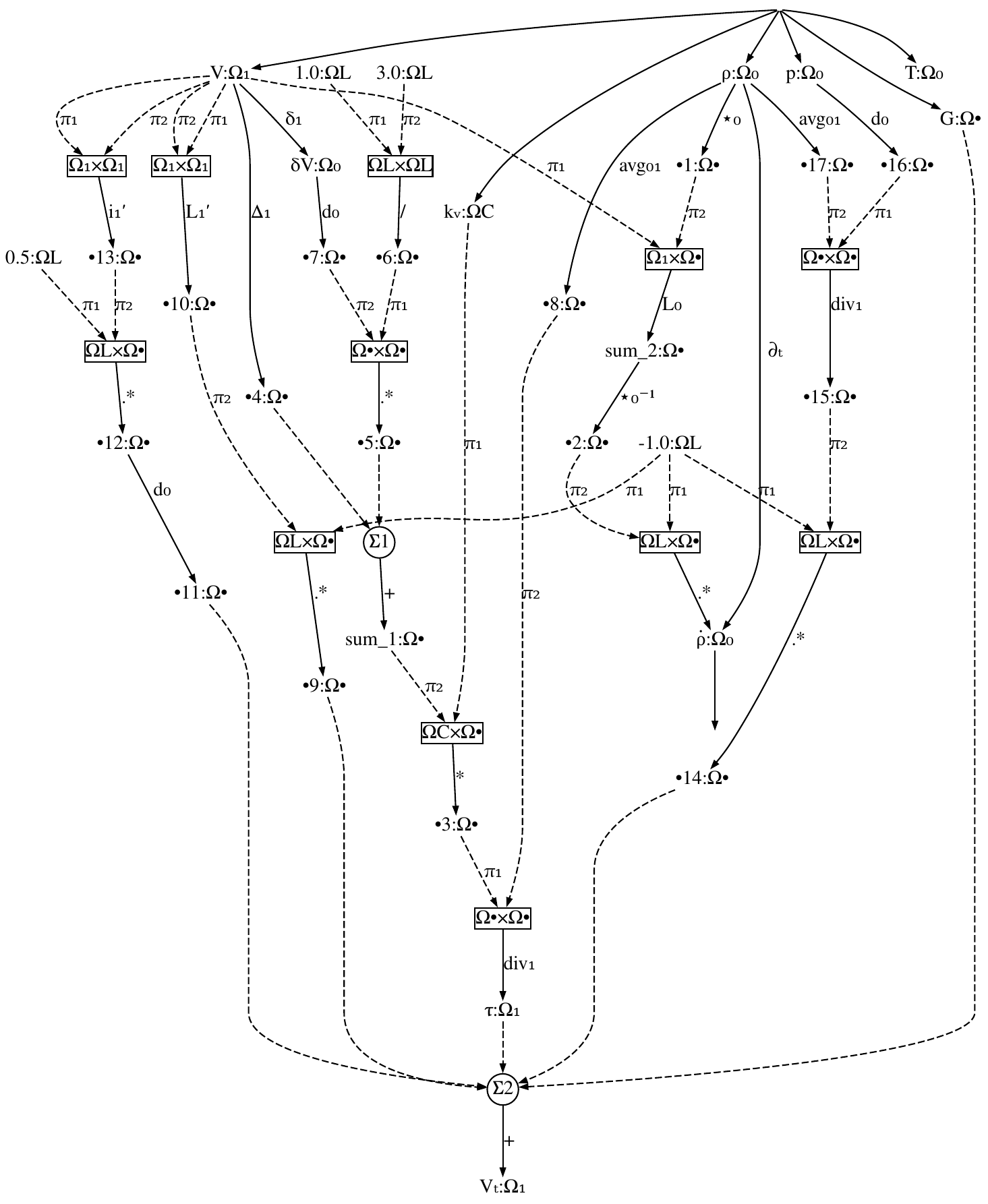}
    \caption{The Buoyancy-Driven Navier-Stokes equations Decapode.}
    \label{fig:bdf_navier_stokes}
    \end{tiny}
\end{figure}

\begin{figure}[htbp]
    \centering
    \begin{subfigure}{0.35\textwidth}
    \includegraphics[width=\textwidth]{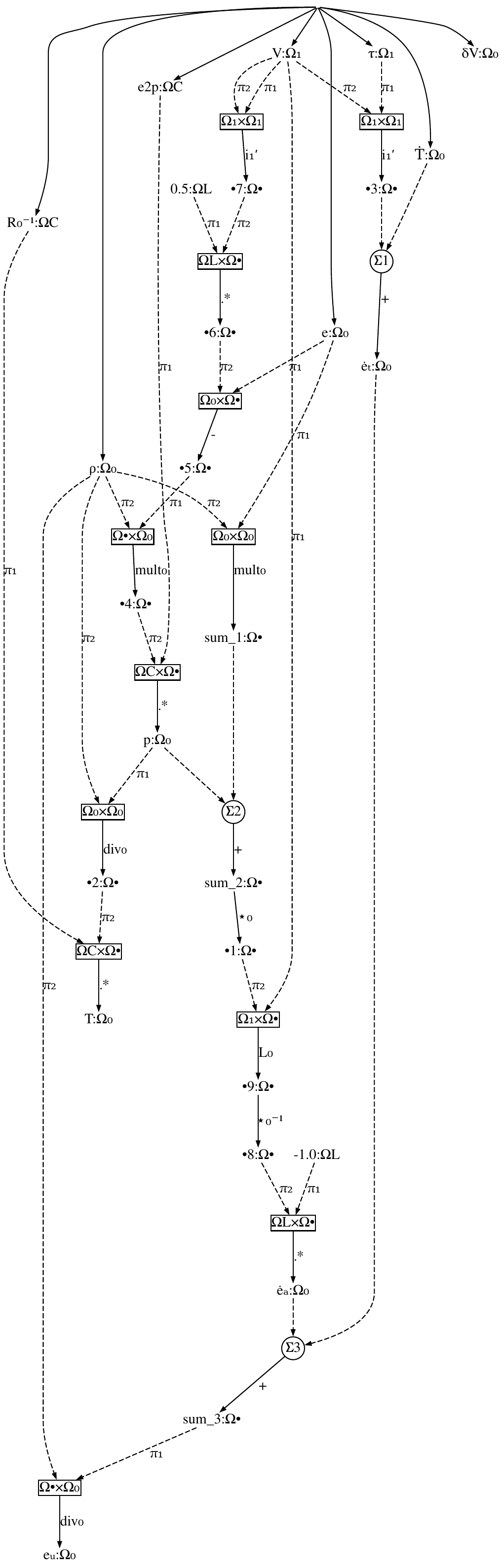}
    \caption{The Buoyancy-Driven energy balance equations Decapode. One may informally observe visually, or formally analyze at the category-theoretic level, the similarities between this formulation and that of Figure ~\ref{fig:CHT_Energy_Decapode}. We see that this formulation is structurally similar, (i.e. in terms of identical sub-Decapodes), but tracks updates in terms of energy, $\mathbf{e}$.}
    \end{subfigure}
    \hspace{0.5cm}
    \begin{subfigure}{0.1\textwidth}
    \includegraphics[width=\textwidth]{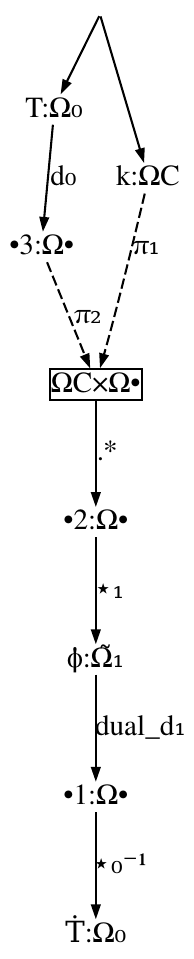}
    \caption{The Buoyancy-Driven flow diffusion equation Decapode.}
    \end{subfigure}
    \hspace{1cm}
    \begin{subfigure}{0.2\textwidth}
    \includegraphics[width=\textwidth]{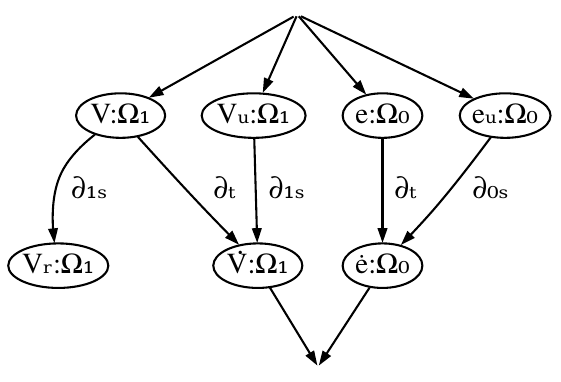}
    \caption{The Buoyancy-Driven flow boundary equations Decapode. Observe here that this component Decapode enforces a boundary condition on the $\mathbf{V}$ primal 1-Form. The composition pattern in Figure ~\ref{fig:cavity_pipe_uwd} is in charge of mapping this bounded version of $\mathbf{V}$ to the values that the Energy and Flow equations internally call $\mathbf{V}$. There is thus a separation of concerns between the components of the composed multiphysics system. The Energy and Flow equations do not need to know whether or when ``boundary logic'' is applied.}
    \label{fig:boundary_eqs}
    \end{subfigure}
    \caption{The energy balance, diffusion, and boundary condition Decapodes for the Buoyancy-Driven fluid flow problem.}
    \label{fig:bdf_supporting_eqs}
\end{figure}

\begin{figure}[htbp]
    \centering
    \includegraphics[width=0.5\textwidth]{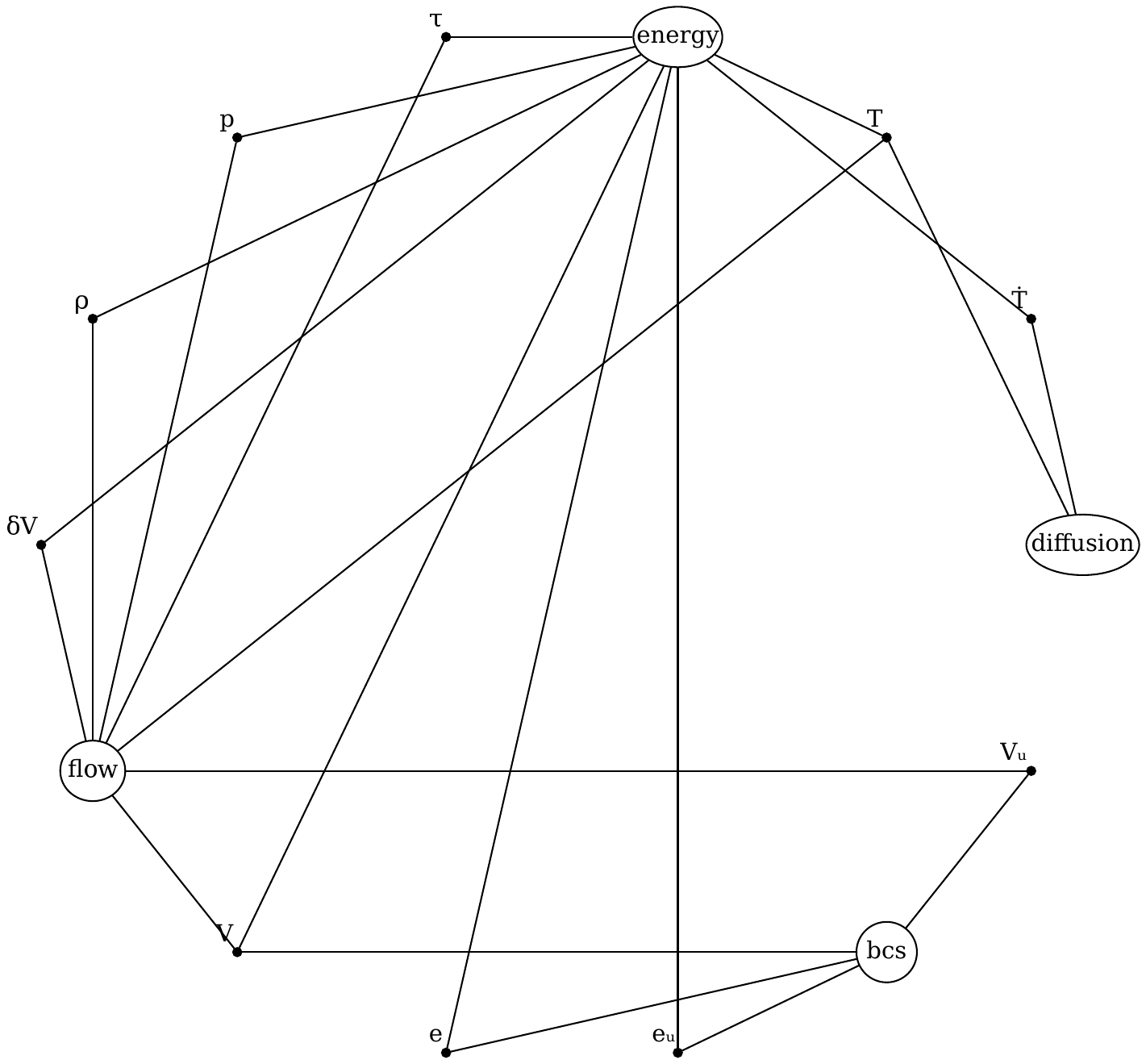}
    \caption{Composition pattern for the Buoyancy-Driven flow benchmark problem. This provides a composition pattern for the Decapodes describing thermal diffusion, the energy balance equations, the fluid flow, and the problem boundary conditions.}
    \label{fig:cavity_pipe_uwd}
\end{figure}

\begin{figure}[htbp]
    \centering
    \begin{tiny}
    \includegraphics[width=0.6\columnwidth]{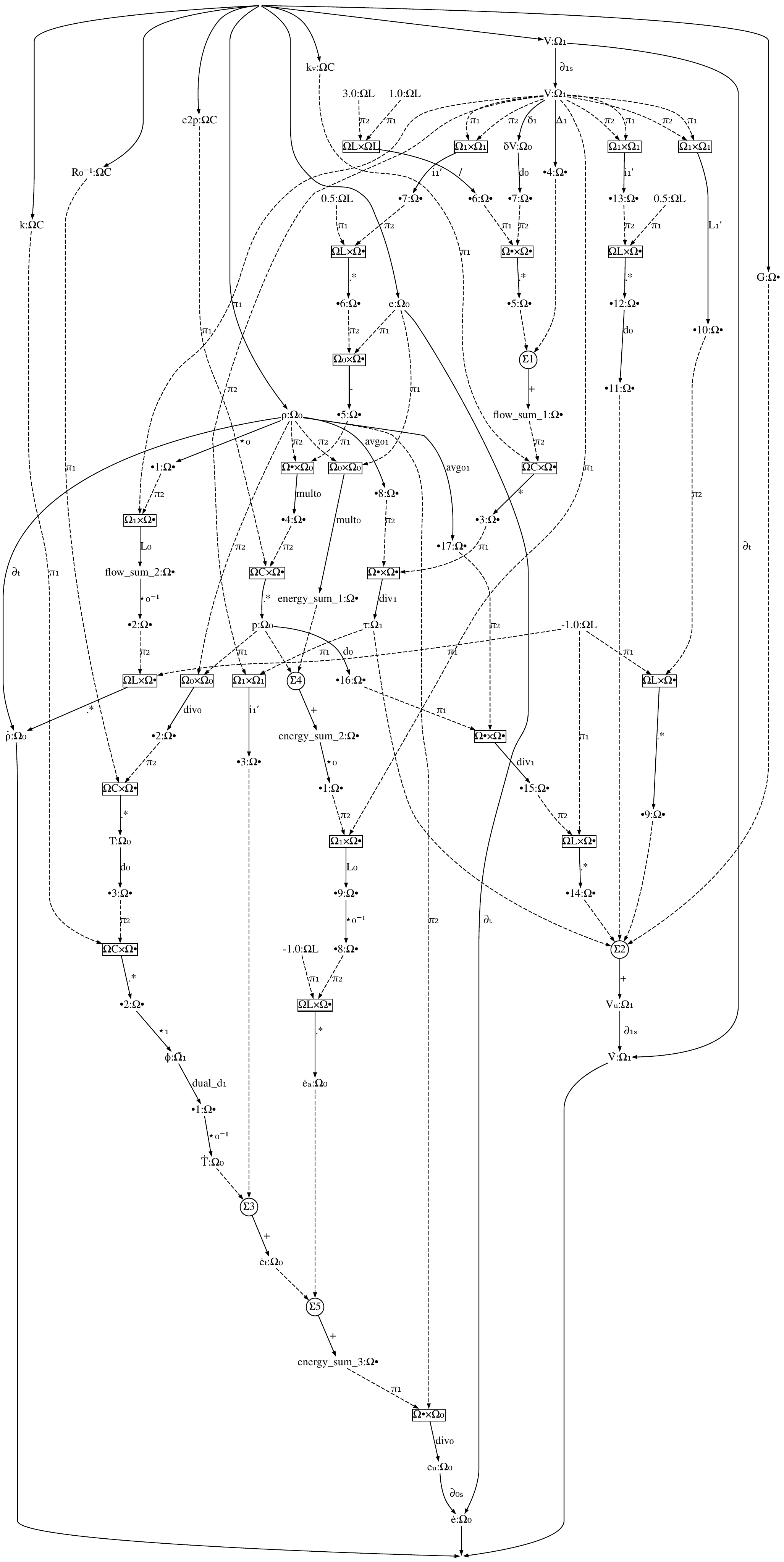}
    \caption{The final composed Buoyancy-Driven flow Decapode. This is an example of a large multiphysics system that would be tedious to write by hand, but by using component Decapodes and applying composition patterns, this complexity is well-managed.}
    \label{fig:cavity_theory}
    \end{tiny}
\end{figure}

As in the previous benchmark, for this problem we use the previously derived Navier-Stokes equation, but this time we still track density, internal energy, and velocity as in Fig.~\ref{fig:bdf_navier_stokes}. We complete this system of equations using the state equations from the SU2 documentation~\footnote{\url{https://su2code.github.io/docs_v7/Theory/\#compressible-navier-stokes}}; Decapodes corresponding to these are given in Fig.~\ref{fig:bdf_supporting_eqs}. The composition pattern of the Decapodes is visualized in  Fig.~\ref{fig:cavity_pipe_uwd}. This results in the Decapode presented in Fig.~\ref{fig:cavity_theory} which is then compiled to a Julia program and run according to the initial conditions and domain described above. The interface conditions are handled by a ``boundary equations Decapode" in Fig.~\ref{fig:boundary_eqs}. These operators apply masks at each discrete step which set the value of $V$ to be zero along the top and bottom walls, enforcing no-slip boundary conditions, and setting the values of $T$ and $p$ to be their initial values, enforcing no change along those walls. The affects of these boundary condition methods on the numerical stability of these simulations is not investigated here, and improving boundary condition application with more advanced techniques is grounds for future work.

\paragraph{Computational Results}

Time marching is performed in SU2 via an implicit solver method to allow for larger timesteps. Again, only explicit solution methods are available for the Decapodes formulation, with $\mathtt{Tsit5}()$ being used from the \texttt{DifferentialEquations.jl} library. The simulations are marched forward for a total time of $t = 5 \mathrm{[s]}$ and the results are compared at steady state. 

\begin{figure}[htbp]
     \centering
     \begin{subfigure}[htbp]{0.49\textwidth}
         \centering
         \includegraphics[width=0.99\textwidth]{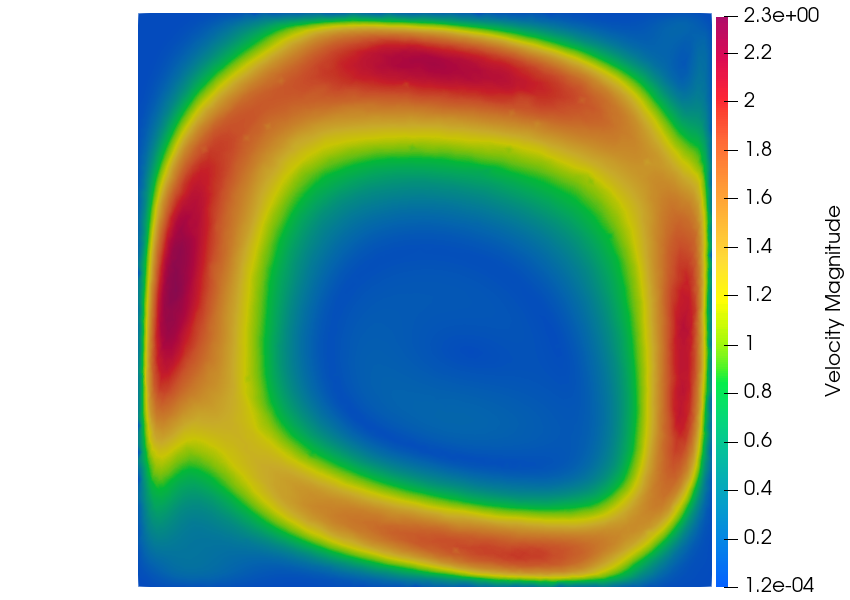}
         \caption{Velocity Field Computed with \Deca}
         \label{fig:vel_dec}
     \end{subfigure}
     \hfill
     \begin{subfigure}[htbp]{0.49\textwidth}
         \centering
         \includegraphics[width=0.99\textwidth]{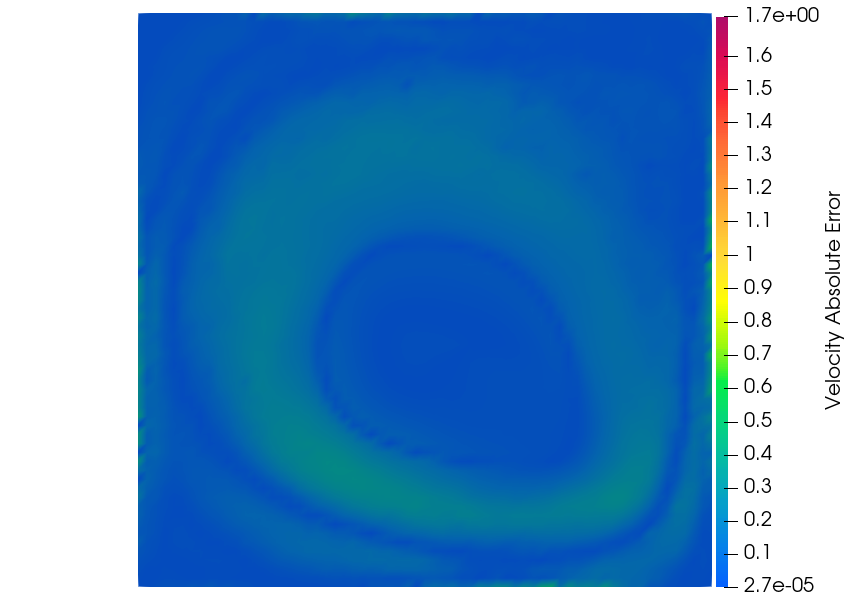}
         \caption{Absolute Error between Velocity Calculation of \Deca and SU2}
         \label{fig:vel_error}
     \end{subfigure}
        \caption{Natural Convection Benchmark Velocity Comparison}
        \label{fig:vel_results}
\end{figure}

\begin{figure}[htbp]
     \centering
     \begin{subfigure}[htbp]{0.49\textwidth}
         \centering
         \includegraphics[width=0.99\textwidth]{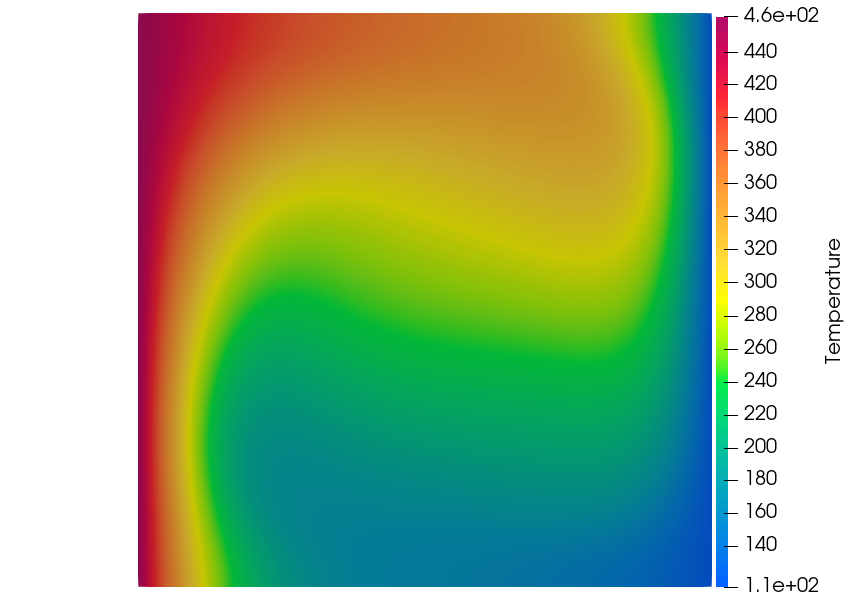}
         \caption{Temperature Field Computed with \Deca}
         \label{fig:temp_dec}
     \end{subfigure}
     \hfill
     \begin{subfigure}[htbp]{0.49\textwidth}
         \centering
         \includegraphics[width=0.99\textwidth]{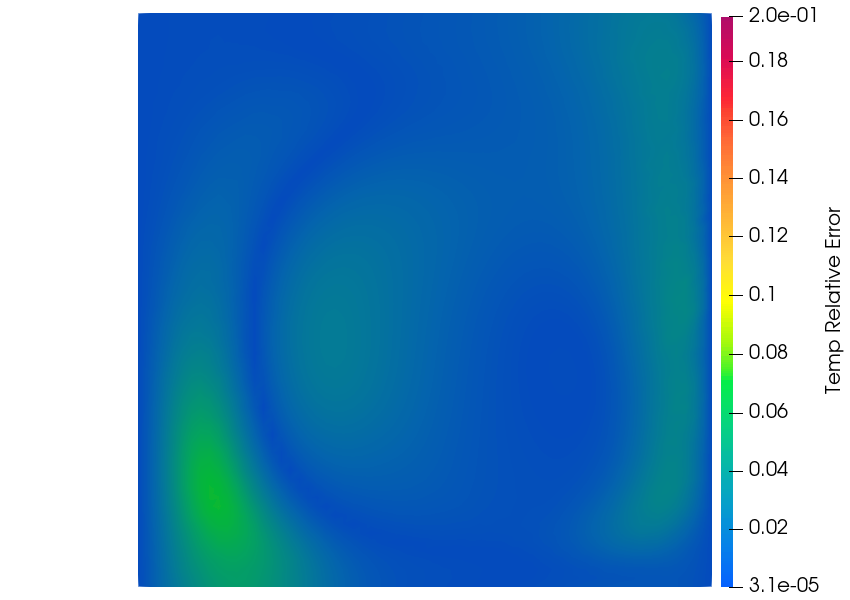}
         \caption{Relative Error between Temperature Calculation of \Deca and SU2}
         \label{fig:temp_error}
     \end{subfigure}
        \caption{Natural Convection Benchmark Temperature Comparison}
        \label{fig:temp_results}
\end{figure}

The qualitative behavior of the velocity and temperature fields are examined in Figures \ref{fig:vel_results} and \ref{fig:temp_results} respectively. Decapodes is generally able to reconstruct the physical behavior in agreement with the results calculated via SU2, with a smaller discrepancy between the two than before. Both tools feature the creation of a clockwork flow and temperature variation in response to bulk changes in density. In Figure \ref{fig:vel_error}, the absolute error between the velocity fields of Decapodes and SU2 is around $\rm 10\%$. Likewise, in Figure \ref{fig:temp_error}, the relative error between the temperature fields calculated by Decapodes and SU2 remains generally low. Absolute and relative error for velocity and temperature can be found in Table~\ref{tab:buoy_temp_benchmarks}.

\begin{figure}[htbp]
     \centering
     \begin{subfigure}[htbp]{0.4\textwidth}
         \centering
         \includegraphics[width=0.99\textwidth]{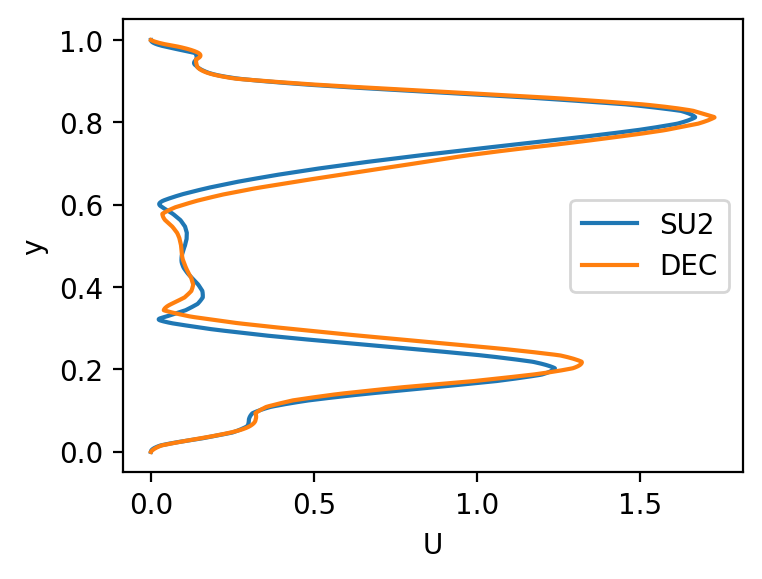}
         \caption{Velocity Profile along midline of domain. Comparison of the velocity magnitude.}
         \label{fig:vel_dec_profile}
     \end{subfigure}
     \hfill
     \begin{subfigure}[htbp]{0.4\textwidth}
         \centering
         \includegraphics[width=0.99\textwidth]{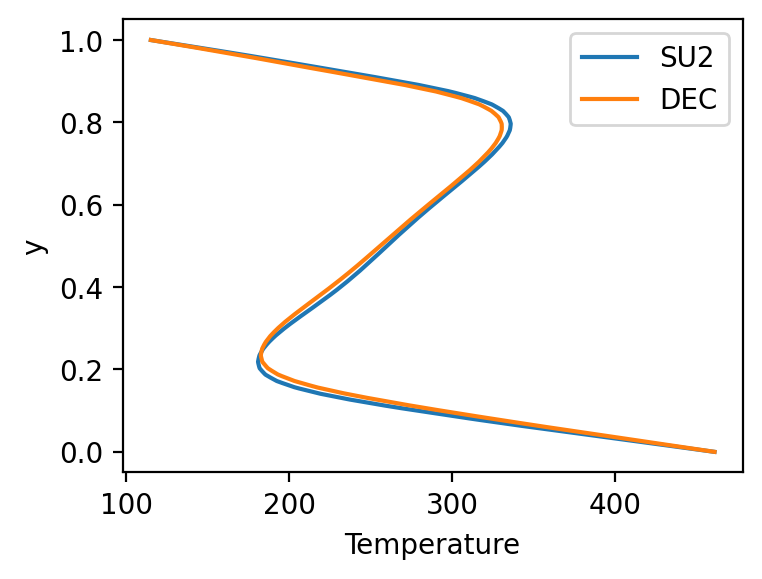}
         \caption{Temperature Profile along midline of domain.}
         \label{fig:vel_error_profile}
     \end{subfigure}
        \caption{Natural Convection Benchmark Velocity Comparison}
        \label{fig:results_profile}
\end{figure}

More granular comparisons between the Decapodes and SU2 results are carried out via vertical profiles at location $x = 0.5$ and presented in Figure \ref{fig:results_profile} for the velocity and temperature profiles. Here, the reconstruction of the behavior is in both qualitative and quantitatively good agreement between the Decapodes and SU2 results, with both having the same general behavior present in the velocity and temperature profiles, and variance of around $10\%$.

\paragraph{Runtime Results}
The SU2 simulation takes $9.08$ seconds of compute time per second of simulation and performs one time step in 
0.0908 seconds. The Decapodes simulation takes 
28600 seconds of compute time per second of simulation and performs one time step in 
0.0721 seconds. As in the previous benchmark, the total difference in computation time is likely dominated by the use of implicit versus explicit time-stepping.

\subsection{Summary of Results}
The same physical phenomena observed in both benchmarks - such as wakes behind cylinders in the CHT benchmark and counter-clockwise flow in the buoyancy-driven flow benchmark - were observed in both the SU2 and Decapodes simulations. Moreover, in both benchmarks, the Decapodes simulations were able to compute results that that agree with those from SU2 up to low average relative errors, as shown in Tables ~\ref{tab:cht_dens_benchmarks} and ~\ref{tab:buoy_temp_benchmarks}. Differences in total computation time are likely dominated by the implicit time-stepping used by SU2 versus explicit time-stepping used by the current Decapodes implementation. We note that there is no property inherent to the Discrete Exterior Calculus that inhibits the use of implicit time-stepping, and so the difference in computation time does not invalidate this approach.

\begin{table}[hbtp]
    \centering
    \begin{tabular}{|l|l|l|}
    \hline
    Experiment & Implementation & Runtime \\
    \hline
    CHT & Decapodes & $50400 \frac{comp}{sim}$ \\[1pt]
    \hline
    CHT & SU2 & $94.3 \frac{comp}{sim}$ \\[1pt]
    \hline
    Buoyancy & Decapodes & $28600 \frac{comp}{sim}$ \\[1pt]
    \hline
    Buoyancy & SU2 & $9.08 \frac{comp}{sim}$ \\[1pt]
    \hline
    \end{tabular}
    \caption{Summary of runtimes for the CHT and Buoyancy experiments in terms of seconds of computation time per second of simulation time. Slower Decapodes times are likely explained by the use of explicit time-stepping.}
    \label{tab:DEC_benchmarks}
\end{table}

\begin{table}[hbtp]
    \centering
    \begin{tabular}{|l|l|l|l|l|l|l|}
    \hline
        Variable & Error & Avg & Max & Min & Total & RMS \\ \hline
        Velocity & Absolute & 0.0991 & 0.560 & 1.83E-5 & 419 & 0.130 \\ \hline
        Temperature & Absolute & 3.50 & 16.7 & 0.000922 & 14800 & 4.69 \\ \hline
        Density & Absolute & 0.00112 & 0.00242 & 0.000604 & 4.75 & 0.00119 \\ \hline
        Velocity & Relative & 0.343 & 32.7 & 0 & 1450 \\ \cline{1-6}
        Temperature & Relative & 0.0139 & 0.0711 & 3.92E-6 & 58.7 \\ \cline{1-6}
        Density & Relative & 0.158 & 0.217 & 0.128 & 667 \\ \cline{1-6}
    \end{tabular}
    \caption{Summary of error results from the Conjugate Heat Transfer experiment.}
    \label{tab:cht_dens_benchmarks}
\end{table}

\begin{table}[!ht]
    \centering
    \begin{tabular}{|l|l|l|l|l|l|l|}
    \hline
        Variable & Error & Avg & Max & Min & Total & RMS \\ \hline
        Velocity & Absolute & 0.0859 & 0.523 & 6.59E-5 & 396 & 0.114 \\ \hline
        Temperature & Absolute & 12.2 & 28.6 & 4.71E-6 & 79200 & 14.6 \\ \hline
        Velocity & Relative & 0.235 & 32.1 & 0 & 1080 \\ \cline{1-6}
        Temperature & Relative & 0.0367 & 0.0861 & 1.63E-8 & 238 \\ \cline{1-6}
    \end{tabular}
    \caption{Summary of error results from the Buoyancy experiment.}
    \label{tab:buoy_temp_benchmarks}
\end{table}

\clearpage

\section{Conclusions}
\label{sec:conclusions}

Decapodes.jl is a software package that combines category theoretic approaches to computer algebra and the discrete differential geometry of the discrete exterior calculus to produce an extremely flexible numerical physics simulation platform. This platform can be used to represent complex multiphysics phenomena and generate performant and accurate numerical simulations of those phenomena. We have demonstrated the utility of this approach with comparisons to the SU2 solver on standard multiphysics problems including conjugate heat transfer and buoyancy-driven flow. These benchmarks show that Decapodes is a viable approach to multiphysics simulation. Additional research is needed to develop more advanced solvers including finite element exterior calculus (FEEC)~\cite{doi:10.1137/1.9781611975543}.

The Decapodes tool represents a paradigm shift in multiphysics modeling. Under this framework, multiphysics simulations can be formally guaranteed to model to the physics they describe and the models themselves are visually-debuggable. Additionally, this tool has drastically reduced the time-to-first-simulation for novel models, enabling a faster iteration loop for physicists. We escape the old ``build-or-buy'' paradigm by automatically generating simulations from a formal description of their physics. By all of this, Decapodes helps democratize the world of multiphysics simulations by lowering the level of effort required for generating complex, accurate simulations.

\appendix

\section{Implementation Choices}
As often is the case for physics simulations, several specific implementation choices were made in the process of creating the simulations shown in Section~\ref{sec:experiments}. Below, we describe the primary implementation choices and custom operators created for these simulations. We believe that as these tools mature and the DEC is developed further, these simulation decisions will both become more rigorous, justified, and result in greater accuracy.

\subsection{Discrete Exterior Calculus}
\label{app:dec_mod}
There were several ways in which the implementation of the DEC diverged from the implementation laid out by Hirani~\cite{hirani2003}. These places of divergence occur for two primary reasons: that Hirani's derivation requires a circumcentric dual and that there were additional operators needed which were not included in Hirani's text. We clarify that this is not due to any lack of clarity or completeness on the part of Hirani, which others have shown is capable of accurate simulation of complex physics~\cite{Mohamed_2016, Mohamed2018NumericalCO, 10.1063/5.0035981}.

As mentioned in Section~\ref{sec:dec}, Hirani's DEC requires a mesh and its dual in order to generate all necessary discrete operators. There are many common ways of generating the dual of a mesh, each with their own benefits and drawbacks. The circumcentric method of mesh subdivision ensures that all dual elements are orthogonal to their corresponding primal elements, an aspect which significantly simplifies several of the discrete operators. Hirani utilizes this method of subdivision because of these mathematical properties, yet in the 2D case this method only provides well-centered duals when the mesh contains no obtuse triangles. The barycentric method of mesh subdivision, on the other hand, will always provide well-centered duals for a 2D triangular mesh.
Although the DEC operators may still of course be instantiated without any guarantee of the mesh being well-centered, through testing, we note that meshes with well-centered duals typically result in more physically-accurate results.
Since the meshes utilized in the examples are irregular and contain obtuse triangles, we used the barycentric subdivision method in our examples. Many of Hirani's operators are dependent on the circumcentric dual assumption, and so we implemented both the geometric hodge star developed in 2021 by Ayoub et. al.~\cite{ayoub_new_2021} and implemented a flat operator which takes the non-orthogonality into account.

Additionally, Hirani's wedge product is defined for products between primal and primal forms, but is not extended to products between primal and dual forms within his thesis. This is instead described in two later papers which describe the special case of a wedge product between a $k$ primal form and $n-k$ dual form where $n$ is the highest simplex dimension and $0 \leq k \leq n$~\cite{desbrun_discrete_2005, seslija_discrete_2011}. While this is not generalized within these papers to wedge products between arbitrary primal and dual forms, it does provide the necessary implementation for calculating inner products.

Finally, we utilize our own implementation of the cross product between primal 1-forms and dual 2-forms in order to calculate self-advection, since that product is necessary for the Lie derivative between two primal 1-forms. There is the potential that the Liebniz rule for the wedge product could provide an implementation for this more general primal-dual wedge product, yet that has not yet been implemented at the time of this writing.

\section*{Acknowledgements}
This work was supported by DARPA under Directly Computable Models: Generalized Algebraic Theories Enhancing Multiphysics Agreement No.~HR00112090067, Young Faculty Award: Model Aware Scientific Computing	Agreement No.~W911NF2010292, ASKEM: Generalized Algebraic Techniques Advancing Scientific Discovery	HR00112220038, and Director's Fellowship: Model Aware Scientific Computing	~W911NF2110323.

\printbibliography
\end{document}